\documentclass[letterpaper,titlepage]{amsart}
\usepackage{amsmath, amssymb, amsthm, amsfonts,xr}
\newtheorem{thm}{Theorem}[section]
\newtheorem{theorem}[thm]{Theorem}

\newtheorem{lemma}[thm]{Lemma}
\newtheorem{lem}[thm]{Lemma}

\newtheorem{corollary}[thm]{Corollary}

\newtheorem{definition}[thm]{Definition}

\newtheorem{problem}[thm]{Problem}
\newtheorem{proposition}[thm]{Proposition}
\newtheorem{remark}[thm]{Remark}

\newtheorem{cor}[thm]{Corollary}

\newtheorem{qu}[thm]{Question}
\newtheorem{prop}[thm]{Proposition}

\newtheorem{dfn}[thm]{Definition}
\newtheorem{rem}[thm]{Remark}

\begin{document}
\title{Dunkl Operators for Complex Reflection Groups}
\author{C.F. Dunkl}
\address{Department of Mathematics\\
University of Virginia\\
Charlottesville\\
Virginia 22904-4137\\
USA}
\email{cfd5z@virginia.edu}
\thanks{During the preparation of this paper, the first-named author
was partially supported by National Science Foundation grant DMS-9970389 and
held a Sesquicentennial Research Associateship at the University of Virginia.}
\author{E.M. Opdam}
\address{Korteweg de Vries Institute for Mathematics\\
University of Amsterdam\\
Plantage Muidergracht 24\\
1018TV Amsterdam\\
The Netherlands}
\thanks{During the preparation of this paper the second named author
was partially supported by a Pionier grant of the Netherlands Organization
for Scientific Research (NWO)}
\email{opdam@science.uva.nl}
\date{June 14, 2001}
\begin{abstract}
Dunkl operators for complex reflection groups are defined
in this paper.
These
commuting operators give rise to a parameter family of deformations of
the polynomial De Rham complex. This leads to the study
of the polynomial ring as a module over the ``rational Cherednik algebra'',
and a natural contravariant form on this module. In the case of the
imprimitive complex reflection groups $G(m,p,N)$, the set of singular
parameters in the parameter family of these structures is
described explicitly, using
the theory of nonsymmetric Jack polynomials.
\end{abstract}
\maketitle
\section{Introduction}
Let $V$ be a finite dimensional Hilbert space.
A finite complex reflection group in the unitary
group $U(V)$ is a finite subgroup which is generated by pseudo-reflections
(or ``complex reflections''),
transformations with all eigenvalues but one equal to unity. In this paper we
introduce for each such group a commutative algebra of parameterized operators
$T_i(k)$ generalizing the partial derivatives $\partial_i=T_i(0)$.
Here $k$ denotes an $N$-tuple of complex parameters where $N$
equals the number of conjugacy classes of complex reflections in the group.

We prove the commutativity by showing that the perturbed De Rham
complex, in which the role of the partial derivatives is played by
the operators $T_i(k)$, is indeed a complex (i.e. the perturbed
differential $d(k)$ satisfies $d(k)^2=0$). At the heart of this argument
lies the computation of a ``Laplacian'' $E(k)$ for this complex.
This operator gives rise to a parameter family of elements in the center of
the group algebra, whose values on the irreducible representations of
$G$ are nonnegative integral linear combinations of the
parameters $k$. These values seem to govern many of the properties of
the operators $T_i(k)$. In particular, analyzing these values easily leads
to a proof of the commutativity of the $T_i(k)$.

In the course of this argument various natural structures arise:
\begin{enumerate}
\item[(i)] The perturbed De Rham complex $P\otimes\bigwedge^\bullet(V^*)$
with differential $d(k)$, where $P$ is the ring of polynomials on $V$.
\item[(ii)] Homogeneous, $G$-equivariant
intertwining operators $\mathcal{V}(k)$ on
$P$, such that $\partial_i\mathcal{V}(k)=
\mathcal{V}(k)T_i(k)$.
\item[(iii)] An hermitian pairing $(\cdot,\cdot)_k$ on $P$ such that
$(x_ip,q)_k=(p,T_i(k)q)_k$.
\item[(iv)] The ring $P$ as a module over the ``rational
Cherednik algebra'' $\mathbb{A}(k)$,
the algebra generated by $\mathbb{C}[G]$, $T_i(k)$ and $P$
(acting on itself by multiplication).
\end{enumerate}
A parameter value $k$ is called {\it singular} if there exists a nonzero
homogeneous polynomial $p$ of positive degree such that for
all $i$, $T_i(k)p=0$.
This turns out to be the only obstruction for the
existence of a homogeneous equivariant intertwining isomorphism as
in (ii). By this remark it is
easy to see that $k$ is singular if and only if there exists an
$i>0$ such that the cohomology group
$H^i(k)$ of the De Rham complex with differential $d(k)$ is nonzero.
Another equivalent formulation is the statement that $(\cdot,\cdot)_k$
is degenerate.
From this one easily sees that $k$ is singular if and only if
$P$ is not irreducible as an $\mathbb{A}(k)$-module.

By these equivalent descriptions it is clear that the set of singular
parameter values is of fundamental importance, and one of the goals
of this paper is to find this set explicitly.
We are not able to solve this problem in general,
but we will derive that the singular set is always a locally finite union
of affine rational hyperplanes in the parameter space.

In the case of a Coxeter group one knows more about the above structures,
and this was described in the paper \cite{DjO}. The present paper grew out
of an attempt to apply the methods discussed in \cite{DjO} to the case
of complex reflection groups.

The complex reflection groups were classified by Shephard and Todd
\cite{ST}. There are 34 exceptional cases called $G_{i}, i=4,..,37$
(containing the exceptional real reflection groups) and an infinite
family of groups $G(m,p,N)$ with $m,p,N\in\mathbb{N}$ and $p|m$.
The group $G\left(m,p,N\right)$ is a subgroup of $U\left(  N\right)$
and consists of
permutation matrices whose nonzero entries are $m^{\mathrm{th}}$ roots of
unity and the product of the nonzero entries is an $\left(  m/p\right)
^{\mathrm{th}}$ root of unity. If $N=1$ we take $p=1$
(cyclic groups of order $m$ acting on $\mathbb{C}$). The groups $G\left(
1,1,N\right) ,\,G\left(  2,1,N\right)  ,\,\allowbreak G\left(  2,2,N\right)
$ are the Coxeter groups of types $A_{N-1},\,B_{N},\,D_{N}$ respectively.

The infinite family $G\left(  m,p,N\right)  $ is studied in detail
in the second half of this paper,
by means of the theory of nonsymmetric Jack polynomials.  A
complete orthogonal decomposition for the pairing associated with $G\left(
m,p,N\right)  $ is obtained with explicit norm formulae. This leads to a
precise description of the set of singular values and the construction of
shift operators which transform between the structures for contiguous
parameter values.

Various alternative interpretations are known for the singular
parameter set in the case of Coxeter groups (see also \cite{DjO}). It
is closely related to the non-semisimple specializations of Hecke algebras
(see \cite{O2}).
Likewise it is closely related to the zeroes of the
Bernstein-Sato polynomial of the discriminant
(see \cite{O1}). In this regard it is also interesting
to compare with the results of the paper \cite{DL}.
It is an interesting question whether some of
these interpretations survive in the general case of
complex reflection groups. The Hecke algebras to be considered
are the topological cyclotomic Hecke algebras studied in
\cite{BMR} (also \cite{O3}).

We note that the algebra $\mathbb{A}(k)$ naturally fits in the
framework of the symplectic reflection algebras that were recently
introduced by Etingof and Ginsburg
\cite{EG}. As is mentioned in that paper, this provides an alternative
proof of the commutativity of the operators $T_i(k)$.
Yet another approach to the proof of the commutativity
is to show the integrability of the related Knizhnik-Zamolodchikov
connection directly, using geometric arguments along the lines of \cite{K}
(see \cite{BMR} for the details of this argument).

\textbf{Acknowledgments.}
The authors wish to thank Pavel Etingof
for a helpful discussion.

The authors were participants in the programme on Symmetric Functions and
Macdonald Polynomials during the preparation of this paper, and
gratefully acknowledge the hospitality of the Isaac Newton Institute
for Mathematical Sciences at Cambridge.

\section{Dunkl operators for complex reflection groups}
\subsection{Complex reflection groups}
Let $V$ be a finite dimensional Hilbert space, and let $U(V)$ be
the group of unitary linear transformations of $V$. An element
$g\in U(V)$ is called a {\it complex reflection} if $g$ has finite
order and $H_g:=\operatorname{Ker}(g-\operatorname{Id})$ is a
hyperplane in $V$. A finite subgroup $G\subset U(V)$ is called a
finite complex reflection group if $G$ is generated by complex
reflections.

Let $G\subset U(V)$ be a finite complex reflection group.
If $W\subset V$ is a
linear subspace, we denote by $G_W\subset G$ the subgroup of those
elements of $G$ which fix the elements of $W$. By a well known result
of Steinberg \cite{St}, $G_W$ is itself a finite complex
reflection group. Clearly, $G_W$ acts faithfully on $W^\perp$. In
particular, if $H\subset V$ is a hyperplane, then $G_H$ is a
cyclic group. If $G_H\not=\{\operatorname{Id}\}$ we call $H$ a
reflection hyperplane. When $H$ is a reflection hyperplane, we
denote by $e_H$ the order of the cyclic group $G_H$. The
collection of all reflection hyperplanes is a central hyperplane
arrangement $\mathcal{A}$ in $V$, on which the group $G$ acts. Let
us denote by $\mathcal{C}$ the set of $G$-orbits in $\mathcal{A}$.
Obviously, $e_H$ only depends on the orbit $C=G\cdot
H\in\mathcal{C}$. We will write $e_C$ instead of $e_H$ whenever
this is convenient.

Let $\operatorname{det}$ be the determinant character on $U(V)$.
When $H\subset V$ is a reflection hyperplane, the characters of
$G_H$ also form a cyclic group, generated by the restriction
$\chi_H$ of $\operatorname{det}$ to $G_H$. We will thus label the
character group of $G_H$ by
\begin{equation}
{\hat{G}}_H=\{\chi_H^{-i}\mid i=0,\dots,e_H-1\}.
\end{equation}
For each of the reflection hyperplanes we choose a functional
$\alpha_H\in V^*$ such that $\operatorname{Ker}(\alpha_H)=H$.

The group $G$ acts on the ring $P$ of polynomials on $V$ in the
usual way, i.e. $p^g(x):=p(g^{-1}x)$. The functional $\alpha_H$
transforms under the action of $G_H$ according to the character
$\chi_H^{-1}$. A polynomial $p$ which transforms under the action of
$G_H$ according to a nontrivial character must vanish on $H$.
Consequently, a polynomial which transforms under the action of
$G_H$ according to $\chi_H^{-i}$ with $i\in\{1,\dots, e_H-1\}$
is divisible by $\alpha_H^i$.

\subsection{Dunkl operators}
Given $H\in\mathcal{A}$ and $i\in\{0,1,\dots, e_H-1\}$, let
\begin{equation}
\epsilon_{H,i}:=\frac{1}{e_H}\sum_{g\in
G_H}\chi_{H}^i(g)g\in\mathbb{C}[G_H],
\end{equation}
be the idempotent of $\mathbb{C}[G_H]$ of the character
$\chi_H^{-i}$. In addition, choose a list of complex numbers $k=(k_{C,i})$,
where $C$ runs over the set of orbits $\mathcal{C}$, and for each
$C\in\mathcal{C}$, $i\in\{1,\dots,e_C-1\}$. With these data we
form, for each reflection hyperplane $H$, the element
\begin{equation}
a_H=a_H(k)=\sum_{i=1}^{e_H-1}e_Hk_{H,i}\epsilon_{H,i}\in\mathbb{C}[G_H].
\end{equation}
Notice that the $\epsilon_{H,i}$ with $i\in\{1,\dots, e_H-1\}$
constitute a basis of the subalgebra in $\mathbb{C}[G_H]$ consisting of
the elements $\sum_{g\in G_H}c_gg\in\mathbb{C}[G_H]$ such that
$\sum c_g=0$. Also notice that the elements $a_H$ are equivariant
with respect to conjugation in the group algebra: for all $g\in
G$,
\begin{equation}
g a_H g^{-1}=a_{gH}.
\end{equation}

Let $\xi\in V$, and denote by $\partial_\xi$ the derivation of
$P$ associated to the constant vector field on $V$ defined by
$\xi$. By what was said above, we can define the following ``Dunkl
operator'' on $P$:
\begin{equation}\label{defdu}
T_\xi(k)=\partial_\xi+\sum_{H\in\mathcal{A}}\alpha_H(\xi)\alpha^{-1}_Ha_H(k),
\end{equation}
where $a_H(k)$ is considered as operator acting on $P$. The operator indeed
maps polynomials to polynomials, since
$\epsilon_{H,0}a_H=0$, so that we can divide by $\alpha_H$
after applying $a_H(k)$. The next proposition is now clear:
\begin{prop}
\begin{enumerate}
\item[(i)] The operator $T_\xi(k)$ does not depend on the choice of
the functionals $\alpha_H$.
\item[(ii)] $T_\xi(k) $ is equivariant with respect to the
action of $G$ on $P$ and $V$:
\begin{equation}
gT_\xi(k)g^{-1}=T_{g\xi}(k).
\end{equation}
\item[(iii)] $T_\xi(k)$ is homogeneous of degree $-1$.
\end{enumerate}
\end{prop}

Let $K^\bullet=P\otimes\bigwedge^\bullet V^*$ denote the algebra of polynomial
differential forms on $V$. Let $\Omega(k)\in \operatorname{End}(P)\otimes K^1$ be given
by
\begin{equation}
\Omega(k)=\sum_{H\in\mathcal{A}}a_H(k)\omega_H,
\end{equation}
where $\omega_H:=\alpha^{-1}_Hd\alpha_H$ is the logarithmic differential of $\alpha_H$ (which
is independent of the choice of $\alpha_H$, and $G_H$-invariant).
\begin{lem}
The operator $\Omega(k):P\to K^1$ is $G$-equivariant with respect to the usual action
of $G$ on $P$, and the diagonal action of $G$ on $K^1$.
\end{lem}
\begin{proof}
This is clear since both $a_H(k)$ and $\omega_H$ are equivariant for the natural actions
of $G$ on $P$ and on the space of (rational) $1$-forms on $V$.
\end{proof}
We write $d(k):P\to K^1$ to denote the map
$d(k)(p)=dp+\Omega(k)(p)$. Thus we have:
\begin{equation}
T_\xi(k)(p)=c_\xi(d(k)(p)),
\end{equation}
where $c_\xi$ denotes the contraction with the constant vector field $\xi$.
We extend the operator $d(k)$ to $K^\bullet$ in the usual way: for all $p\in P$
and $\omega\in \bigwedge^\bullet V^*$ we define $d(k)(p\otimes\omega)=d(k)(p)\wedge\omega$,
and extend this linearly to $K^\bullet$. Note however that, unlike the case $k=0$,
this is not a derivation of the algebra $K^\bullet$.
\begin{lem}
We have the following equivalent definition of $d(k)$ on $K^\bullet$.
Extend $\Omega(k)$ to the $G$-equivariant endomorphism
$\Omega(k)$ on $K^\bullet$  which is defined by $\Omega(k)(\omega)\allowbreak =
\sum_{H\in\mathcal{A}}a_H(k)(\omega_H\wedge\omega)$ for
$\omega\in K^\bullet$, where $a_H(k)$
acts diagonally on $K^\bullet$.
Then $d(k)=d+\Omega(k)$.
The operator $d(k)$ is equivariant for the
diagonal action of $G$ on $K^\bullet$.
\end{lem}
\begin{proof}
If $g\in G_H$, then $x^g\in x+\mathbb{C}\alpha_H$, for all $x\in V^*$.
In addition, $\omega_H$ is $G_H$-invariant.
Hence if $\omega=p\otimes dx_1\wedge\dots\wedge dx_l$, then
\begin{equation}
(a_H(k)(p))\otimes\omega_H\wedge dx_1\wedge\dots\wedge dx_l=
a_H(k)(p\otimes\omega_H\wedge dx_1\wedge\dots\wedge dx_l),
\end{equation}
where we used the diagonal action of $G$ on the right hand side.
The equivariance of $d(k)$ follows from the previous Lemma.
\end{proof}

Consider the Koszul differential $\partial$ on $K^\bullet$. This differential is
defined by:
\begin{align*}
  \partial:K^l &\to K^{l-1} \\
  p\otimes dx_{1}\wedge\dots\wedge dx_{l} &\to \sum_{r=1}^l (-1)^{r+1}x_{r}p\otimes
  dx_{1}\wedge\dots\wedge\widehat{dx_{r}}\wedge\dots\wedge dx_{l}.
\end{align*}
Observe that $\partial$ is $U(V)$-equivariant with
respect to the diagonal action of $U(V)$ on $K^\bullet$.

Let $E(0)$ denote the Euler vector field on $V$. This vector field is the infinitesimal
generator of the action of $\mathbb{C}^\times$ on $V$ (by scalar multiplication).
Differentiating the diagonal action of $\mathbb{C}^\times$ on $K^\bullet$ we obtain
an action of $E(0)$ on $K^\bullet$, the ``diagonal action''. In other words,
if we put $K^l_m=P_m\otimes \bigwedge^lV^*$, then
$E(0)$ has eigenvalue $l+m$ on
$K^l_m$. Notice that $d(k)(K^l_m)\subset K^{l+1}_{m-1}$ and that
$\partial(K^l_m)\subset K^{l-1}_{m+1}$.
\begin{prop}
Let $E(k)=E(0)+\sum_{H\in\mathcal{A}}a_H(k)$, acting
diagonally on $K^\bullet$.
Then $\partial d(k)+d(k)\partial=E(k)$.
\end{prop}
\begin{proof}
It is well known that $\partial d(0)+d(0)\partial=E(0)$,
since $d(0)$ is the ordinary
De Rham differential on $K^\bullet$.

Extend the operators $d(k)$ and $\partial$ to the complex
$\overline{K}^\bullet$ of rational differential
forms on $V$ in the natural way. For each $H\in \mathcal{A}$, define the operator
$w_H:\overline{K}^\bullet\to\overline{K}^\bullet$ by $w_H(\eta):=\omega_H\wedge\eta$.
By the previous Lemma and the equivariance of $\partial$ we have
\begin{align*}
&\partial a_H w_H+a_H w_H\partial=\\
&a_H(\partial w_H+w_H\partial)=a_H.
\end{align*}
This finishes the proof.
\end{proof}
The next Lemma is of crucial importance in all that follows.
\begin{lem}\label{lem:real}
We put $z(k)=\sum_{H\in\mathcal{A}}a_H(k)\in\mathbb{C}[G]$.
This element of the group algebra has
the following properties:
\begin{enumerate}
\item[(i)] The element $z(k)$ is in the center of $\mathbb{C}[G]$.
\item[(ii)] For $(V_\tau,\tau)\in \hat{G}$, let $c_\tau(k)$ denote the scalar such that $z(k)$
acts on $V_\tau$ by multiplication with $c_\tau(k)$. Then $c_\tau(k)$ is a linear function of
$k$, with nonnegative integer coefficients.
\item[(iii)] Let $\operatorname{triv}$ denote the trivial representation of $G$. Then
$c_\tau(k)\equiv 0$ if and only if $\tau=\operatorname{triv}$.
\end{enumerate}
\end{lem}
\begin{proof}
(i) This follows immediately from the equivariance of the elements $a_H$.

(ii) Let the restriction of $\tau$ to $G_H$ be
\begin{equation}
\tau|_{G_H}=\sum_{i=0}^{e_H-1}n_{H,i}^\tau\chi_H^{-i}
\end{equation}
for certain nonnegative integers $n_{H,i}^\tau$. Observe that these branching numbers only depend
on the orbit $C=G\cdot H$ of $H$. Hence the trace of $z(k)$ on $V_\tau$ equals
\begin{equation}
\operatorname{trace}\tau(z(k))=\sum_{C\in\mathcal{C}}
\sum_{i=1}^{e_H-1}|C|e_Cn_{C,i}^\tau k_{C,i},
\end{equation}
showing that $c_\tau(k)$ is linear with nonnegative rational numbers as coefficients.
On the other hand, $z(k)$ is a linear expression in the $k_{C,i}$, with
coefficients that are central elements of the algebra $A[G]$, where $A$ denotes the
ring of algebraic integers. By a well known result in the theory of representations of
finite groups,
the central elements of $A[G]$ assume algebraic integer values on the irreducible representations
of $G$ over $\mathbb{C}$. This proves the result.

(iii) If $c_\tau(k)\equiv 0$ then it follows from the proof of (ii) that
for each $H\in\mathcal{A}$, $\tau|_{G_H}$ contains only the trivial
character of $G_H$. Since $G$ is generated by the subgroups $G_H$, this implies that $\tau=
\operatorname{triv}$.
\end{proof}
\begin{cor}
Let $X,Y,Z$ denote the generators of an associative algebra $\tilde{C}$, satisfying the
relations $XY+YX=Z$, $X^2=[Z,X]=[Z,Y]=0$. The map $X\to \partial$, $Y\to d(k)$
and $Z\to E(k)$ extends to a representation of $\tilde{C}$ on $K^\bullet$.
\end{cor}
\begin{proof}
If $\tau\in\hat{G}$, we put $K^l_{m,\tau}$ for the $\tau$-isotypic component of $K^l_m$.
We have $\partial(K_{m,\tau}^l)\subset K_{m+1,\tau}^{l-1}$,
$d(k)(K_{m,\tau}^l)\subset K_{m-1,\tau}^{l+1}$, and finally
$E(k)(K_{m,\tau}^l)\subset K_{m,\tau}^l$. It also follows immediately that $E(k)$ commutes with $\partial$
and $d(k)$.
\end{proof}
We put $K(r,\tau):=\oplus_{l+m=r}K^l_{m,\tau}$, and $K(0)=K^0_0=\mathbb{C}$.
\begin{cor}
$K(r,\tau)$ is a finite dimensional $\tilde{C}$-submodule, and $K^\bullet$ decomposes as direct sum
$K^\bullet=\oplus_{r\geq 0,\tau}K(r,\tau)$.
\end{cor}

Furthermore we write $K^\bullet(+):=\oplus_{r>0,\tau}K(r,\tau)$, which is the $\tilde{C}$-submodule
of $K^\bullet$ complementary to $K(0)$.
We thus have the decomposition $K^\bullet=K^\bullet(+)\oplus K(0)$.
\begin{cor}\label{cor:ker}
Assume that for all $\tau\in \hat{G}$, $-c_\tau(k)\not\in\mathbb{N}$.
Then
\begin{equation}
\operatorname{Ker}(d(k))\cap K^\bullet(+)=\operatorname{Ker}(d(k))\cap\operatorname{Im}(d(k)).
\end{equation}
In particular, $\operatorname{Ker}(d(k))\cap K^0(+)=\{0\}$.
\end{cor}
\begin{proof}
$E(k)$ acts by scalar multiplication with $r+c_\tau(k)$ on the submodule $K(r,\tau)$.
Therefore, by the above assumption, $E(k)$ is invertible on each $\tilde{C}$-submodule $K(r,\tau)$,
with the exception of $K(0)$. Write $E(k)^{-1}$ for the inverse of $E(k)$ on $K^\bullet(+)$.
Let $\omega\in\operatorname{Ker}(d(k))\cap K^\bullet(+)$, then
\begin{align*}
\omega&=E(k)^{-1}E(k)(\omega)\\
&=E(k)^{-1}(\partial d(k)+d(k)\partial)(\omega)\\
&=d(k)(\partial(E(k)^{-1}\omega)).\\
\end{align*}
This shows that $\eta=\partial(E(k)^{-1}\omega)\in K^\bullet(+)$ is a solution of
$d(k)(\eta)=\omega$, proving the desired result.
\end{proof}
\begin{thm}\label{def}
Assume that for all $\tau\in \hat{G}$, $-c_\tau(k)\not\in\mathbb{N}$.
We call a linear operator $\mathcal{V}:K^\bullet\to K^\bullet$ completely homogeneous if
$\mathcal{V}(K^l_m)\subset K^l_m$, for all $l,m$.
There exists a unique completely homogeneous linear map $\mathcal{V}(k):K^\bullet\to K^\bullet$ such that
\begin{enumerate}
\item[(i)] $\mathcal{V}(k)$ is the identity operator on $K(0)=\mathbb{C}$,
\item[(ii)] $\mathcal{V}(k)(p\otimes\omega)=(\mathcal{V}(k)p)\otimes\omega$, for all
$\omega\in\bigwedge^\bullet V^*=K_0^\bullet$, and
\item[(iii)] $d(k)\mathcal{V}(k)=\mathcal{V}(k)d(0)$.
\end{enumerate}
Moreover, $\mathcal{V}(k)$ is a $G$-equivariant linear isomorphism.
\end{thm}
\begin{proof}
First we show that if $\mathcal{V}(k)$ exists, it is necessarily a linear isomorphism.
If not, let $m$ be minimal such that $\mathcal{V}(k)$ has a nontrivial kernel in $K_m^\bullet$. By
(ii) this implies that $P_m\cap\operatorname{Ker}(\mathcal{V}(k))\not=0$, and by (i) we see that $m>0$.
Let $0\not=p\in P_m$ be such that $\mathcal{V}(k)p=0$. By (iii) we see $\mathcal{V}(k)d(0)p=0$.
By the assumption on $m$ this
implies that $d(0)p=0$, and since $m>0$ this is a contradiction. Note that this
argument is independent of the assumption on $k$, since it only uses that $d(0)p=0$ implies that $p=0$
if $p\in P_m$ with $m>0$.

A similar argument shows that $\mathcal{V}(k)$ must be unique.
If not, there exists a nonzero completely homogeneous operator
$\mathcal{W}$ satisfying (ii) and (iii) but with $\mathcal{W}(K_0^\bullet)=0$.
Let $m>0$ be minimal such that $\mathcal{W}(K_m^\bullet)\not=0$.
By (ii) this implies that $\mathcal{W}(P_m)\not=0$. Let $p\in P_m$ such that
$\mathcal{W}p\not=0$. But then
$d(k)\mathcal{W}p=\mathcal{W}d(0)p=0$, which implies that
$\mathcal{W}p\in\operatorname{Ker}(d(k))\cap K^0(+)$.
Given the assumption on $k$, this contradicts Corollary \ref{cor:ker}.

We now construct $\mathcal{V}(k)$ by induction on the degree $m$.
Suppose that $m>0$ and that we have already constructed
$\mathcal{V}(k)$ on $K^\bullet_i$ for $i<m$, satisfying (i), (ii) and (iii). Let $p\in K^0_m=P_m$.
Then $d(k)(\mathcal{V}(k)d(0)p)=0$, and thus, by the previous corollary, there exists a unique
$q\in P_m$ such that $d(k)q=\mathcal{V}(k)d(0)p$. Hence we define $\mathcal{V}(k)p=q$. For any $\omega\in K_0$
and $p\in P_m$ we now put $\mathcal{V}(k)(p\otimes\omega)=(\mathcal{V}(k)p)\otimes\omega$, and use this to define
$\mathcal{V}(k)$ on $K^\bullet_m$. It is immediate that $\mathcal{V}(k)$ satisfies (i), (ii) and (iii).

The $G$-equivariance of $\mathcal{V}(k)$ now follows from the equivariance of $d(0)$ and $d(k)$.
The equivariance implies that
for any $g\in G$, $g\circ \mathcal{V}(k)\circ g^{-1}$ also meets the requirements (i), (ii) and (iii). By the
uniqueness property we conclude that $\mathcal{V}(k)=g\circ \mathcal{V}(k)\circ g^{-1}$.
\end{proof}
\begin{cor}\label{acy}
The map $d(k)$ is a differential on $K^\bullet$, i.e. $d(k)^2=0$.
In addition, the complex $(K^\bullet(+),d(k))$ (with
$K^\bullet(+):=\oplus_{r>0,\tau}K(r,\tau)$) is acyclic (i.e. its cohomology is $0$) if
we assume that for all $\tau\in \hat{G}$, $-c_\tau(k)\not\in\mathbb{N}$.
\end{cor}
\begin{cor}\label{clif}
The representation of $\tilde{C}$ on the submodule $K(r,\tau)$ (in which $Z=E(k)$ acts by scalar
multiplication with $s:=r+c_\tau(k)$) factors through the Clifford algebra $C(s)$, the
associative algebra with generators $X$ and $Y$, and
relations $XY+YX=s$, $X^2=Y^2=0$.
\end{cor}
The next reformulation is the main result of this section:
\begin{thm}
For all $\xi,\eta\in V$,  $T_\xi(k)T_\eta(k)=T_\eta(k)T_\xi(k)$.
\end{thm}
\begin{proof}
Choose a basis $e_i$ of $V$, with dual basis $x_j$ of $V^*$. Put $T_i=T_{e_i}(k)$.
A simple direct computation shows that for all $p\in P$,
\begin{equation}
d(k)^2p=\sum_{i<j}(T_iT_j-T_jT_i)p\otimes dx_i\wedge dx_j.
\end{equation}
Hence the statement $d(k)^2=0$ is equivalent to the commutativity of the $T_i$.
\end{proof}
\begin{cor}
The restriction of $\mathcal{V}(k)$ to $P$ satisfies $T_\xi(k)\mathcal{V}(k)=\mathcal{V}(k)\partial_\xi$.
\end{cor}
\begin{proof}
The point is that, because of defining property (ii) of $\mathcal{V}(k)$, $\mathcal{V}(k)$ commutes with
the contraction $c_\xi$ for each $\xi\in V$. Thus we obtain:
\begin{align*}
T_\xi(k)\mathcal{V}(k)p&=c_\xi(d(k)\mathcal{V}(k)p)\\
&=c_\xi(\mathcal{V}(k)d(0)p)=\mathcal{V}(k)\partial_\xi p.\\
\end{align*}
\end{proof}

When $s:=r+c_\tau(k)\not=0$, the Clifford algebra $C(s)$ is semisimple and has
only one irreducible module $M^s\simeq\mathbb{C}^2$, with basis $a,b$ such
that
\begin{equation}
X=\begin{pmatrix}
    0 & 0 \\
    1 & 0
  \end{pmatrix},
Y=
  \begin{pmatrix}
    0 & s \\
    0 & 0
  \end{pmatrix}.
\end{equation}
Thus $K(r,\tau)$ is isomorphic to a direct sum of copies of the irreducible representation
$\tau\otimes M^s$ of $\mathbb{C}[G]\otimes C(s)$. The ``cohomology'' $H(\tau\otimes M^s):=
\operatorname{Ker}(Y)/\operatorname{Im}(Y)$ is equal to $0$ in the module $\tau\otimes M^s$.
However, the algebra $C(0)$ is the Grassmann algebra. This algebra is no longer semisimple.
The space of equivalence classes of indecomposable modules of the Grassmann algebra
is complicated, and the ``cohomology'' $\operatorname{Ker}(Y)/\operatorname{Im}(Y)$
is not necessarily $0$ in the indecomposable modules. This happens for the only
irreducible representation of $C(0)$, the trivial representation, but it may happen for
nontrivial indecomposable modules as well. Note that for all parameters
values $k$, the trivial representation of $C(0)$ is contained in $K^\bullet$ at least once,
in the form of the submodule $K(0)$.
\subsection{Singular parameter values}\label{sub:sing}
Let $H^i(k)$ denote the i-th cohomology group of the complex
$(K^\bullet,d(k))$.
The following are equivalent:
\begin{cor}\label{cor:coh}
\begin{enumerate}
\item[(i)] $H^i(k)=0$, for all $i>0$.
\item[(ii)] $H^0(k)=\mathbb{C}$.
\item[(iii)] There exists a completely homogeneous intertwining map
$\mathcal{V}(k)$ satisfying the conditions (i), (ii) and (iii) of Theorem \ref{def}.
\end{enumerate}
\end{cor}
\begin{proof}
(i)$\Rightarrow$(ii): The complex $(K^\bullet, d(k))$ is a direct sum of the
finite dimensional subcomplexes $K^\bullet(n)=\oplus_{l+m=n}K^l_m$. Let us
denote by $H^i(n,k)$ the cohomology groups of $(K^\bullet(n),d(k))$.
The Euler characteristic
\begin{equation}
\chi(n):=\sum_{i\geq 0}(-1)^i\operatorname{dim}(H^i(n,k))=
\sum_{i\geq 0}(-1)^i\operatorname{dim}(K^i_{n-i})
\end{equation}
is independent of $k$. Hence for all $k$, $\chi(n)=0$ if $n>0$.
Thus for all $n>0$ and $k$ we have: if $H^i(n,k)=0$ for all $i>0$,
then $H^0(n,k)=0$.

(ii)$\Rightarrow$(iii): We carry out the proof of the existence of an intertwining map
$\mathcal{W}(k)$ as in
Theorem \ref{def}, but with the role of $d(0)$ and $d(k)$ interchanged. We can do
this, because we now know that $d(k)^2=0$. Since we assume that $H^0(k)=\mathbb{C}$,
we see that the proof of Theorem \ref{def} showing that the intertwiner has to be a linear
isomorphism also applies in this situation.
Hence $\mathcal{V}(k)=\mathcal{W}(k)^{-1}$ is an intertwining operator as required.

(iii)$\Rightarrow$(i): The argument in the proof of Theorem \ref{def} showing that $\mathcal{V}(k)$
has to be a linear isomorphism applies for all $k$.
Hence $\mathcal{V}(k)$ defines an isomorphism between
the cohomology spaces $H^i(0)$ and $H^i(k)$.
\end{proof}
\begin{dfn}
The parameter $k$ is called regular if the equivalent statements of
\ref{cor:coh} hold.
Otherwise we say that $k$ is singular.
\end{dfn}
We saw in Corollary \ref{acy} that if $k$ is singular,
then there exists a $\tau\in\hat{G}$ such that
$-c_\tau(k)=m\in\mathbb{N}$. More precisely we have:
\begin{cor}\label{cor:hyp} For $\tau\in\hat{G}$, let $P_\tau$ denote the
$\tau$-isotypical component of $P$, and $P_{m,\tau}$ the $\tau$-isotypical
component of $P_m$.
If $k$ is singular, there exists a $\tau\in\hat{G}$ and $m\in\mathbb{N}$
such that $c_\tau(k)+m=0$ and $P_{m,\tau}\not=0$.
\end{cor}
\begin{proof}
If $k$ is singular then there exists a submodule $K(m,\tau)$ of $K^\bullet$
such that
$H^0(K(m,\tau),d(k))\not=0$ and $m\in\mathbb{N}$.
This implies that $(P_m:\tau)>0$ and $c_\tau(k)+m=0$.
\end{proof}
The converse statement is false.
A counterexample occurs already in the case of the symmetric
group $S_4$. Let $\tau$ be the irreducible 3 dimensional
representation $(2,1,1)$. Then $P_{5,\tau}\not=0$ and $c_\tau(k)=8k$,
but $k=-5/8$ is not a singular parameter.

In the case of a finite Coxeter group, the set of singular $k$ can be
determined
exactly, cf. \cite{DjO}.
In the present generality we do not know how to prove a similar explicit
description.
\subsection{An hermitian form}
For $\xi\in V$, denote by $\xi^*\in V^*$ the
element such that for all $\eta\in V$, $\xi^*(\eta)=(\xi,\eta)$. The map $\xi\to\xi^*$ is an anti-linear
isometric isomorphism. We extend this to an anti-linear isomorphism $*:S\to P$, where $S$
denotes the symmetric algebra on $V$. The inverse map is denoted by $*$ as well.
By the commutativity of the operators $T_\xi(k)$, we can uniquely extend the linear map
$\xi\to T_\xi(k)$ to obtain a linear map from $S$ to $\operatorname{End}(V)$. This map will
be denoted by $s\to s(T)$ (or by $s\to s(T(k))$ if necessary).
We define a sesquilinear pairing $(\cdot, \cdot)_k$ on $P$ by
\begin{equation}\label{eq:pair}
(p,q)_k:=(p^*(T(k))q)(0).
\end{equation}
\begin{prop}\label{prop:easy}
\begin{enumerate}
\item[(i)] For all $g\in G$, $(p^g, q^g)_k=(p,q)_k$.
\item[(ii)] $(\xi^*p,q)_k=(p,T_{\xi}(k)q)_k$.
\item[(iii)] $(P_{m,\tau}, P_{l,\sigma})_k=0$ if $m\not=l$ or if $\tau\not=\sigma$.
\end{enumerate}
\end{prop}
\begin{proof}
(i) This follows from the remark that for all $p\in P$, $g\in G$, $(p^g)^*=g(p^*)$, since
for we have, for all $x\in V^*$ and $v\in V$:
\begin{equation}
(gx^*,v)=(x^*, g^{-1}v)=x(g^{-1}v)=x^g(v)=((x^g)^*,v).
\end{equation}
Hence by the equivariance of $T_\xi(k)$ we have:
\begin{align*}
(p^g,q^g)_k&=((g(p^*)(T(k)))q^g)(0)\\
&=((g\circ p^*(T(k))\circ g^{-1})q^g)(0)=(p^*(T(k))q)^g(0)=(p,q)_k.
\end{align*}

(ii) This follows immediately from the definition, using the commutativity of the $T_\xi(k)$.

(iii) From the definition we see that $(P_m, P_l)_k=0$ if $l\not=m$.
By (i) it follows that $(ep,q)_k=(ep,eq)_k=(p,eq)_k$ for every
self-adjoint idempotent $e$ of $\mathbb{C}[G]$.
This implies (iii).
\end{proof}
The main theorem of this section is:
\begin{thm}\label{thm:herm}
The pairing satisfies $(p,q)_k={\overline{(q,p)_{\overline{k}}}}$.
\end{thm}
\begin{proof}
We prove this by induction on the degree $m$. By Proposition \ref{prop:easy} it is
enough to show that for every
$\tau\in\hat{G}$ and $m\in \mathbb{Z}_+$,
$(p,q)_k={\overline{(q,p)_{\overline{k}}}}$ if $p,q\in P_{m,\tau}$.
Assume by induction that this statement holds true for all $p,q\in P_{l,\sigma}$ with $l<m$
(the case $m=0$ being trivial). Let $e_i$ denote an orthonormal basis of $V$,
and let $x_i=e_i^*$ denote the dual basis of coordinates on $V$.
Then, since $\overline{c_\tau(k)}=c_\tau(\overline{k})$ by
Lemma \ref{lem:real}, we have for all $p,q\in P_{m,\tau}$:
\begin{align*}
(m+c_\tau(k))(p,q)_k&=(E({\overline{k}})p,q)_k\\
&=(\sum_ix_iT_i({\overline{k}})p,q)_k\\
&=\sum_i(T_i({\overline{k}})p,T_i(k)q)_k\\
&=\sum_i{\overline{(T_i(k)q,T_i({\overline{k}})p)_{\overline{k}}}}\\
&=\sum_i{\overline{(x_iT_i(k)q,p)_{\overline{k}}}}\\
&=(m+c_\tau(k)){\overline{(q,p)_{\overline{k}}}}.
\end{align*}
Since $(p,q)_k$ depends polynomially on $k$ and since $(l+c_\tau(k))\not=0$
for generic
values of $k$, this proves the necessary induction step.
\end{proof}
\begin{cor}
For all $x\in V^*$ and $p,q\in P$
we have $(T_{x^*}({\overline{k}})p,q)_k=(p,xq)_k$.
\end{cor}
\begin{proof}
We have
\begin{align*}
(T_{x^*}({\overline{k}})p,q)_k&=
{\overline{(q,T_{x^*}({\overline{k}})p)_{\overline{k}}}}\\
&={\overline{(xq,p)_{\overline{k}}}}\\
&=(p,xq)_k.
\end{align*}
\end{proof}
As was noticed in Proposition \ref{prop:easy},
the finite dimensional subspaces $P_{m,\tau}\subset P$ satisfy
$(P_{m,\tau},P_{l,\sigma})_k=0$ unless $m=l$ and $\tau=\sigma$.
It follows that there exists
a polynomial $p\not=0$ such that $(p,P)_k=0$ if and only if
there exists a polynomial $q\not=0$ such that $(P,q)_k=0$.
In this case we call the sesquilinear pairing $(\cdot,\cdot)_k$
degenerate.
\begin{prop}\label{prop:equiv}
The following are equivalent:
\begin{enumerate}
\item[(i)] $k$ is singular.
\item[(ii)] $(\cdot, \cdot)_k$ is degenerate.
\item[(iii)] There exists a proper graded ideal $I\subset P$ which is stable
for the action of the operators $T_\xi(k)$.
\item[(iv)] ${\overline{k}}$ is singular.
\end{enumerate}
\end{prop}
\begin{proof}
(i)$\Rightarrow$(ii) Choose $m>0$ such that there exists a $0\not=p\in P_m$ with $d(k)p=0$.
Then $(P,p)_k=0$.

(ii)$\Rightarrow$(iii) Take $I:=\{p\mid (P,p)_k=0\}$.

(iii)$\Rightarrow$(i) Let $m>0$ be minimal such that $I_m\not=0$. Then $d(k)I_m=0$.

(i)$\Leftrightarrow$(iv) By the above text, there exists a $q\not=0$
such that $(P,q)_k=0$ if and only if there exists a $p\in P$, $p\not=0$
such that $(p,P)_k=0$. But by Theorem \ref{thm:herm} this is also equivalent
to $(P,p)_{\overline{k}}=0$. Hence, using the equivalence of (i) and (ii),
we see that (i) is indeed equivalent to (iv).
\end{proof}
\begin{cor}
The set $K^{sing}$ of singular parameter values consists of an
infinite union of hyperplanes $K_{m,\tau}$. Here $K_{m,\tau}$ denotes
the hyperplane given by the equation
$m+c_\tau(k)=0$, and the union runs only over pairs $(m,\tau)$ with
$\tau$ nontrivial and $P_{m,\tau}\not=0$.
\end{cor}
\begin{proof}
By Proposition \ref{cor:hyp} we know that $K^{sing}$ is contained
in the above union of hyperplanes $K_{m,\tau}$.
Notice that this is a locally finite set of hyperplanes.
On the other hand, by the above Proposition,
$K^{sing}=\cup_{m\in\mathbb{N}}K^{degen}_m$ where $K^{degen}_m$ denotes
the set of parameter values $k$ such that $(\cdot,\cdot)_k$ is degenerate
on $P_m$. Now
$K^{degen}_m$ is given by the condition that the determinant of
$(\cdot,\cdot)_k$ with respect to a basis of $P_m$ vanishes. Hence
$K^{degen}_m$ is an algebraic hypersurface in the parameter space.
Since it has to be contained in $\cup K_{m,\tau}$, it follows
that, as a subset of the parameter space, $K^{degen}_m$ is a
union of hyperplanes $K_{m,\tau}$. This proves the claim.
\end{proof}
\begin{cor}
Let $K_r$ denote the $\mathbb{Q}$-vector space of rational
parameters $k$
(i.e. $k_{C,i}\in\mathbb{Q}$ for all $C$, $i$).
The set $K^{sing}\cap K_r$
is a (locally finite) union of hyperplanes, and
$K^{sing}$ is the union of the complexifications of these rational
hyperplanes.
\end{cor}
\begin{proof}
The hyperplanes $K_{m,\tau}$ are all rational, by Lemma \ref{lem:real}.
\end{proof}
\begin{cor}
For every pair $(m,\tau)$, there exists a basis $(b_i)$
of $P_{m,\tau}$ such that the determinant of $(b_i,b_j)_k$
is a product of linear factors of the form $l+c_\sigma(k)$
with $\sigma\in \hat{G}$ nontrivial, and $P_{l,\sigma}\not=0$.
\end{cor}
\begin{proof}
The determinant is a polynomial in $\mathbb{C}[K]$ such that
its zero set is contained in $K^{degen}_m$. Hence it contains
only irreducible factors of the form $l+c_\sigma(k)$. The determinant
is determined up to multiplication by an arbitrary nonzero positive
real number by change of basis. When we specialize at $k=0$, the pairing
is clearly positive definite hermitian, and thus we can fix the
normalization by choosing an appropriate basis.
\end{proof}
The above results show that, in order to describe the
set $K^{sing}$, it suffices to describe $K^{sing}\cap K_r$.
By Proposition \ref{prop:equiv} this is equal to the
set $K^{degen}\cap K_r$. In particular we may restrict to real
parameters, which has the advantage that,
by Theorem \ref{thm:herm},
the form $(\cdot,\cdot)_k$ is hermitian:
\begin{prop}\label{prop:pos}
Suppose that $k$ is real. Then
\begin{enumerate}
\item[(i)]
$(\cdot, \cdot)_k$ is hermitian.
\item[(ii)] Suppose moreover
that, for all nontrivial irreducible representations
$\tau\in\hat{G}$, $c_\tau(k)+m(\tau)>0$, where
$m(\tau)$ denotes the lowest homogeneous degree
such that $(\tau:P_{m(\tau)})\not=0$
(this holds in particular when all the parameters satisfy $k_{C,j}\geq 0$).
Then $(\cdot,\cdot)_k$ is positive definite.
\end{enumerate}
\end{prop}
\begin{proof}
(i) This is immediate from Theorem \ref{thm:herm}.

(ii)
This follows by induction on the homogeneous degree $m$, by taking $p=q$ in
the computation in the proof of Theorem \ref{thm:herm}.
\end{proof}
\subsection{Lowest weight modules over the rational Cherednik algebra}
We assume that $k$ is real throughout this subsection.
Let us consider the structure of $P$ as a module over the rational
Cherednik algebra $\mathbb{A}(k)$,
the algebra generated by $\mathbb{C}[G]$, $T_\xi(k)$ and
$P$ (acting on itself by multiplication) (see \cite{EG}).
\begin{lem} For $p\in P$, denote by $m(p)$ the operator
$m(p):P\to P$, $m(p)(q)=pq$ (multiplication by $p$).
For convenience, we put $k_{C,0}=0$, and for $j\in\mathbb{Z}$ we define
$k_{C,j}=k_{C,j^\prime}$ if $j-j^\prime$ is divisible by $e_C$.
In $\mathbb{A}(k)$ we have
\begin{equation}\label{eq:com}
[T_\xi(k),m(p)]=m({\partial_\xi p})+
\sum_{H\in\mathcal{A}}\sum_{i,j=0}^{e_H-1}e_H(k_{H,i+j}-k_{H,j})
\alpha_H(\xi)\alpha_H^{-1}
m({\epsilon_{H,i}(p)})\epsilon_{H,j}.
\end{equation}
In particular, we have
\begin{equation}\label{eq:def}
[T_\xi(k),m(x)]=x(\xi)+\sum_{H\in\mathcal{A}}\sum_{g\in G_H,g\not= 1}
c_g(k)\alpha_H(\xi)x(v_H)\alpha_H(v_H)^{-1} g,
\end{equation}
where $v_H\in V$ is the vector such that $(v_H,v)=\alpha_H(v)$ for all $v\in V$,
and where $c_g(k)$ is the constant
\begin{equation}\label{eq:cg}
c_g(k)=\sum_{j=0}^{e_H-1}\chi^j_H(g)(k_{H,j+1}-k_{H,j}).
\end{equation}
The function $g\to c_g(k)$ is invariant for conjugation of $g$ by $G$.
\end{lem}
\begin{proof}
This follows in a straightforward way from the equations
\begin{equation}
p=\sum_{i=0}^{e_H-1}\epsilon_{H,i}(p)
\end{equation}
and
\begin{equation}
[\epsilon_{H,j},m({\epsilon_{H,i}(p)})]=m({\epsilon_{H,i}(p)})
(\epsilon_{H,j-i}-\epsilon_{H,j}).
\end{equation}
\end{proof}
Let us describe how $\mathbb{A}(k)$ fits into the framework of the paper \cite{EG}.
Consider the abstract associative algebra $T(V+V^*)\otimes\mathbb{C}[G]$,
the smash product where $G$ acts diagonally on the tensor algebra $T(V+V^*)$.
We introduce in this algebra the relations $[\xi,\eta]=0$ (for all $\xi, \eta
\in V$), $[x,y]=0$ (for all $x,y\in V^*$) and finally,
\begin{equation}
[\xi,x]=x(\xi)+\sum_{H\in\mathcal{A}}\sum_{g\in G_H,g\not= 1}
c_g(k)\alpha_H(\xi)x(v_H)\alpha_H(v_H)^{-1} g,
\end{equation}
where $c_g(k)$ is as in equation \ref{eq:cg}.
The resulting algebra $\mathbb{A}^\prime(k)$ is a symplectic reflection algebra
in the sense of \cite{EG}.
In particular, by Theorem 1.3 of \cite{EG}, $\mathbb{A}^\prime(k)$ has the PBW-property
(this means that $\mathbb{A}^\prime(k)$ is isomorphic as a vector space to
$P\otimes S\otimes \mathbb{C}[G]$). By construction, $\mathbb{A}(k)$ is a quotient
of $\mathbb{A}^\prime(k)$ via $P\ni p\to m(p)$, $S\ni p^*\to p^*(T(k))$ and
$g\to g$ (action in $P$). Using the PBW-property one easily identifies the
$\mathbb{A}^\prime(k)$-module $P$ as
\begin{equation}
P=\operatorname{Ind}_{S\otimes\mathbb{C}[G]}^{\mathbb{A}^\prime(k)}(triv),
\end{equation}
where $triv$ is the one dimensional representation such that $G$ acts trivially, and
$triv(V)=0$.
In fact, it is not hard to see that $P$ is a faithful module over $\mathbb{A}^\prime(k)$
(see \cite{EG}, Proposition 4.5). We will therefore identify $\mathbb{A}^\prime(k)$
and $\mathbb{A}(k)$ from now on.

In the above construction we identified $P$ with an induced module. We generalize this
construction in the following way.
Let $(V,\tau)$ be an irreducible module of $G$. We extend $\tau$ to the algebra
$S\otimes \mathbb{C}[G]$ by demanding that $\tau(V)=0$. We define
\begin{equation}
M(\tau,k):=\operatorname{Ind}_{S\otimes\mathbb{C}[G]}^{\mathbb{A}^\prime(k)}(\tau).
\end{equation}
In the sequel we will usually suppress the parameter $k\in K$ in the notation if
no confusion is possible.
As a vector space, $M(\tau)\simeq P\otimes V$. The action of $G$ is the diagonal
action, $P$ acts by multiplication in the left factor of the tensor product, and
the action of $T_\xi(k)$ is given by
\begin{equation}\label{eq:act}
T_\xi(p\otimes v)=
{\partial_\xi p}\otimes v+
\sum_{H\in\mathcal{A}}\sum_{i,j=0}^{e_H-1}e_H(k_{H,i+j}-k_{H,j})
\alpha_H(\xi)\alpha_H^{-1}
{\epsilon_{H,i}(p)}\otimes\epsilon_{H,j}(v),
\end{equation}
according to equation \ref{eq:com}. So we have in particular that $P=M(triv)$.
\begin{lem}\label{lem:eigenv}
Let $\sigma\in\hat{G}$ and $m\in \mathbb{Z}_+$.
Let $M(\tau)_{\sigma}$ of denote the $\sigma$-isotypic component
of $M(\tau)$, and let $M(\tau)_{m,\sigma}=M(\tau)_{\sigma}\cap (P_m\otimes V)$.
The deformed Euler vector field $E(k)=\sum_i x_iT_i(k)$ acts
on $M(\tau)_{m,\sigma}$ by multiplication with the scalar
$m-c_\tau(k)+c_\sigma(k)$.
\end{lem}
\begin{proof}
From equation \ref{eq:act} we
see that, for $p\in P_m$,
\begin{align*}
E(k)(p\otimes v)
&=mp\otimes v
-\sum_{H\in \mathcal{A}}
\sum_{j=0}^{e_H-1}e_Hk_{H,j}p\otimes \epsilon_{H,j}(v)\\
&{\ \ \ \ \ \ \ \ \ \ \ \ \ \,}+\sum_{H\in\mathcal{A}}\sum_{l=0}^{e_H-1}
e_Hk_{H,l}\sum_{i+j=l\operatorname{mod}e_H}
\epsilon_{H,i}(p)\otimes\epsilon_{H,j}(v)\\
&=(m-c_\tau(k))p\otimes v+\sum_{l=0}^{e_H-1} e_Hk_{H,l}
\epsilon_{H,l}^\Delta(p\otimes v),\\
\end{align*}
where $\epsilon^\Delta_{H,l}$ denotes the idempotent $\epsilon_{H,l}$ acting
diagonally on $P\otimes V$. Thus on the $\sigma$-isotypical part
$M(\tau)_{m,\sigma}$ of $P_m\otimes V$,
the action of $E(k)$ is scalar with eigenvalue
$m-c_\tau(k)+c_\sigma(k)$, as claimed.
\end{proof}
\begin{prop}\label{prop:red}
All $\mathbb{A}(k)$-submodules of $M(\tau)$ are graded. In other words,
if $M\subset M(\tau)$ is an $\mathbb{A}(k)$-submodule, then
\begin{equation}
M=\oplus_{m\in \mathbb{Z}_+} M_m
\end{equation}
with $M_m:=M\cap (P_m\otimes V)$. With respect to this grading,
$T_\xi(k)$ has degree $-1$, $x$ has degree $+1$, and $g\in G$ has degree $0$.
\end{prop}
\begin{proof}
Let $M\subset M(\tau)$ be an $\mathbb{A}(k)$-submodule.
Since $\mathbb{C}[G]\subset\mathbb{A}(k)$,
we have $M=\oplus_{\sigma\in\hat{G}}M_\sigma$.
By the preceding Lemma, the eigenvalues of the
operator $E(k)$ separate elements of different homogeneous degree
in each isotypic part $M_\sigma$.
Therefore $M_\sigma=\oplus_nM_\sigma\cap(P_m\otimes V)$, and thus $M$ itself
is also the direct sum of its graded pieces.
\end{proof}
\begin{cor} The module $M(\tau)$ has a unique proper maximal submodule.
In particular, $M(\tau)$ is indecomposable.
\end{cor}
\begin{proof} Since $\tau$ is irreducible for the action of $G$,
a submodule $M\subset M(\tau)$ is proper if and only if $M\cap V=\{0\}$
(where $V=M(\tau)_0$ as before).
All submodules $M\subset P$ are graded, and therefore the sum $M^\prime$ of all
proper submodules of $M(\tau)$ also has the property that $M^\prime\cap V=\{0\}$.
Hence $M^\prime$ is the unique maximal proper submodule of $M(\tau)$.
\end{proof}
\begin{dfn}
We call a module $M$ over $\mathbb{A}(k)$ a
lowest weight module with lowest $G$-type $\tau$
if $M$ is a nontrivial quotient of $M(\tau,k)$. We denote
by $L(\tau,k)$ (or simply $L(\tau)$)
the unique irreducible quotient of $M(\tau,k)$.
\end{dfn}
The above Proposition \ref{prop:red} shows that all lowest weight
modules $M$ with lowest $G$-type $\tau$ have a unique ``natural'' grading
\begin{equation}
M=\oplus_{n\in\mathbb{Z}_+}M_n
\end{equation}
such that $M_0$ is the irreducible $\mathbb{C}[G]$-module of type $\tau$. Note
that the submodules of lowest weight modules are also graded (as they are
subquotients of $M(\tau)$).
Note however that the grading induced by $M$ on a lowest weight submodule
$M^\prime$ of $M$ is shifted with respect to the natural grading of
$M^\prime$ (unless $M=M^\prime$).

Notice the analogy with the theory of highest weight modules for a semisimple Lie
algebra over $\mathbb{C}$. The role of Verma-modules is played by the modules
$M(\tau)$.
Let us call an element $m\in M$ of a lowest weight module $M$ of $\mathbb{A}(k)$
{\it primitive} if $T_\xi(k)m=0$ for all $\xi\in V$.
Clearly the subspace $M^p$ of primitive
elements in $M$ is a graded $G$-subspace of $M$.
\begin{prop}\label{prop:hom}
Let $M$ be an $\mathbb{A}(k)$-module generated by a subspace $M_0$ of
primitive vectors such that $M_0$ is an irreducible $\mathbb{C}[G]$-module
of type $\tau$. Then there exists a surjective homomorphism $\phi:M(\tau)\to M$ such
that $\phi(M(\tau)_0)=M_0$. The homomorphism $\phi$ is unique up to a scalar.
\end{prop}
\begin{proof}
This is clear by the universal property of induced modules.
\end{proof}
\begin{prop}\label{prop:emb}
\begin{enumerate}
\item[(i)] Let $M$ be a lowest weight module with lowest $G$-type $\tau$, graded with
its natural grading. Then the $\sigma$-isotypic component $M^p_{\sigma}$ of $M^p$
is contained
in $M_{m(\tau,\sigma)}$ with $m(\tau,\sigma):=c_\tau(k)-c_\sigma(k)$.
In particular, $M^p$ is finite dimensional.
\item[(ii)] Let $m\in\mathbb{Z}_+$ be such that $M^p_{m}\not=0$, and
let $H\subset M^p_m$ be an irreducible $G$-subspace of $M^p_m$ of type $\sigma$.
The subspace
$J(H):=PH\subset M$ is a lowest weight submodule of $M$, of lowest $G$-type
$\sigma$.
\item[(iii)] Let $m\in\mathbb{Z}_+$ be maximal such that $M^p_{m}\not=0$, and
let $H\subset M^p_m$ be an irreducible $G$-subspace of $M^p_m$ of type $\sigma$.
Then $J(H)\simeq L(\sigma)$.
\item[(iv)]
If, for all $\sigma\in\hat{G}$, $c_\tau(k)-c_\sigma(k)\not\in \mathbb{N}$,
then $M(\tau)$ is irreducible.
\end{enumerate}
\end{prop}
\begin{proof}
(i) Suppose $M^p_{m,\sigma}\not=0$.
By definition of primitivity, $E(k)(M^p_{m,\sigma})=0$. On the other
hand, $E(k)$ acts by multiplication with the scalar $m-c_\tau(k)+c_\sigma(k)$ on
$M(\tau)_{m,\sigma}$.
Since $M_{m,\sigma}$ is a quotient of this space, $E(k)$ also acts on $M_{m,\sigma}$
via multiplication with this scalar. The equation $m=m(\tau,\sigma)$
for the degree $m$ follows. In particular, the dimension of
$M^p$ is bounded by
\begin{equation}\label{eq:dm}
d(M):=\sum_{\sigma\in\hat{G}}\operatorname{dim}(M_{m(\tau,\sigma)}).
\end{equation}

(ii) $J(H)$ is an $\mathbb{A}(k)$-submodule of $M$, because
$T_\xi(k)(ph)=[T_\xi(k),m(p)](h)$ and, by equation \ref{eq:com},
we have $[T_\xi(k),m(p)]\in P\otimes\mathbb{C}[G]$. By the above Proposition,
$J(H)$ is a quotient of $M(\sigma)$.

(iii) By the condition on $m$ we see that
$J(H)^p=H$, since clearly $J(H)^p\subset M^p$. Thus every nonzero submodule of
$J(H)$ contains $H$, and is therefore equal to $J(H)$. Thus $J(H)$ is the unique
irreducible quotient of $M(\sigma)$.

(iv) This is a special case of (iii), since the condition implies that
the maximal value of $m$ such that $M(\tau)^p_m\not=0$ is equal to $0$.
Thus we can take $H=V=M(\tau)_0$ in (ii) to see that $M(\tau)=L(\tau)$
in this case.
\end{proof}
\begin{cor}
Each lowest weight module $M$ has a finite Jordan-H\" older series whose irreducible
quotients are isomorphic to modules of the form $L(\sigma)$.
\end{cor}
\begin{proof}
By proposition \ref{prop:emb} it is clear that $M$ contains an irreducible
submodule $J_1$ isomorphic to $L(\sigma_1)$ for some $\sigma_1$. If
$N_1:=M/J_1=0$ we are done. If not, we continue by choosing
an irreducible submodule $J_2\simeq L(\sigma_2)$ in
$N_1$, and form the quotient $N_2:=N_1/J_2$. We thus construct a sequence of consecutive
quotients $M\to N_1\to N_2\to\dots$. Notice that by construction,
$d(N_i)<d(N_{i-1})$ (see equation \ref{eq:dm}) in each step of the process.
Hence the process has to stop in finitely many steps, say $N_n=0$. Now put
$M_i:=\operatorname{Ker}(M\to N_{n-i})$. Then we get
\begin{equation}
M=M_0\supset M_1\supset M_2\supset\dots\supset M_n=0,
\end{equation}
with $M_i/M_{i+1}=\operatorname{Ker}(N_{n-i-1}\to N_{n-i})\simeq L(\sigma_{n-i})$,
as desired.
\end{proof}
We thus arrive at the following fundamental
\begin{problem}\label{pbl} For a lowest weight module $M$ over $\mathbb{A}(k)$,
denote by $[M:L(\sigma,k)]$ the
multiplicity of $L(\sigma)$ in a Jordan-H\" older series of $M$. Compute
and interpret the multiplicities $[M(\tau,k):L(\sigma,k)]$.
\end{problem}
There are many natural questions and problems related to the topics introduced in this
subsection. One should study the structures that were introduced before in the case of
$P=M(triv)$ (such as the contravariant pairing, the De Rham complex, the singular set)
systematically for general lowest weight modules. In addition, one should consider
the natural analogue of category $\mathcal{O}$ in the present situation.
Except for the remarks below, we resist the temptation to address
any of these questions here, since this would take us too far afield.
The main goal of this paper is the study
of the special case $\tau=triv$.

Some straightforward remarks are in order. By Proposition \ref{prop:emb} it is clear
that if $M$ has lowest $G$-type $\tau$, $[M:L(\sigma,k)]$ is bounded by
the multiplicity $[M_{m(\tau,\sigma)}:\sigma]$ of $\sigma$ in the degree
$m(\tau,\sigma)=c_\tau(k)-c_\sigma(k)$ part of $M$. In particular we see that
\begin{equation}
[M(\tau,k):L(\tau,k)]=1
\end{equation}
and that for $\sigma\not=\tau$,
\begin{equation}
[M(\tau,k):L(\sigma,k)]\not= 0\text{\ implies\ that\ }c_\tau(k)-c_\sigma(k)\in\mathbb{N}.
\end{equation}
Hence if we introduce an ordering of $\hat{G}$ by defining $\sigma\geq\tau$ if
and only if $c_\sigma(k)\leq c_\tau(k)$, then the matrix $[M(\tau):L(\sigma)]$ is
unipotent upper triangular. In particular, the matrix is invertible.

Define the $\tau$-singular set
$K^{sing}_\tau\subset K$ as the set of $k\in K$ such that $M(\tau,k)\not=L(\tau,k)$.
The description of $K_\tau^{sing}$ is just one part of Problem \ref{pbl}.
By the above it is plain that $k\not\in K^{sing}_\tau$ if for all $\sigma$,
$c_\tau(k)-c_\sigma(k)\not\in\mathbb{N}$. The converse is in general not true (as we
saw in the special case $\tau=triv$). We see that $K_\tau^{sing}$ is contained
in a locally finite collection of hyperplanes, such that the coefficients of
the affine linear forms describing the hyperplanes are integral (but with signs).
In the case where $\tau$ is the linear character of $G$
whose restriction to $G_H$ is $\chi_H^{-b_\tau(C)}$ (for some $b_\tau(C)
\in\{0,\dots, e_H-1\}$) if $H\in C$, then it follows by formula \ref{eq:act} that
$M(\tau,k)\simeq M(triv,\beta_\tau(k))$ where ${\beta_\tau(k)}_{C,j}:=
k_{C,j+b_\tau(C)}-k_{C,b_\tau(C)}$. In particular, $K^{sing}_\tau=
\beta^{-1}_\tau(K^{sing})$.
\subsection{The module $P$ over $\mathbb{A}(k)$}
Let us return to the module $P=M(triv)$.
In \cite{DjO} a polynomial $q\in P_+$ was called {\it singular} for the parameter
$k$ if $T_\xi(k)q=0$ for all $\xi\in V$. In other words, $q$ is singular if and only if
$q\in H^0(k)$ (see subsection \ref{sub:sing}) and $q(0)=0$.
In the language of the previous subsection, a polynomial is singular if and only if it is a
primitive element of positive degree for the module $P$. In addition, $H^0(k)=P^p$.

Recall that, by definition, nonzero singular polynomials for $k$
exist if and only if $k$ is a singular parameter.
The space $H^0(k)$ of primitive polynomials  for $k$ is graded and is a $G$-space.
For each $\tau\in \hat{G}$ and degree $m\in\mathbb{Z_+}$, we denote by $H^0_{m,\tau}(k)$
the space of primitive polynomials for $k$ in degree $m$ and of type $\tau$.
By Proposition \ref{prop:emb}, we have the
relation
\begin{equation}
m=-c_\tau(k)
\end{equation}
if $H^0_{m,\tau}\not=0$,
showing that $H^0(k)$ is finite dimensional.
Suppose that $H\subset H^0_{m,\tau}(k)$ is an irreducible subspace.
By Proposition \ref{prop:hom}, the ideal $J(H):=PH$
it generates in $P$ is a lowest weight $\mathbb{A}(k)$-submodule of $P$ with lowest $G$-type
$\tau$.

The form $(\cdot,\cdot)_k$ on $P$
is an hermitian contravariant form on the module $P$ over $\mathbb{A}(k)$, in
the sense of Proposition \ref{prop:easy} (i) and (ii). This construction can be extended
to all lowest weight modules, and it plays
a role quite similar to the Shapovalov form on Verma-modules for a semisimple Lie
algebra. We will restrict ourselves to the special case at hand, the module $P$.
The radical $\operatorname{Rad}(k):=\{p\in P\mid (p,P)_k=0\}$ of this
form is a graded ideal of $P$,
which is stable for $G$ and for the application of
the operators $T_\xi(k)$ (see Proposition \ref{prop:equiv} and its proof).
In other words, it is a (graded) $\mathbb{A}(k)$-submodule of $P$.
\begin{prop}
The radical $\operatorname{Rad}(k)$ of $(\cdot,\cdot)_k$
is the unique maximal proper $\mathbb{A}(k)$-submodule of $P$.
In particular, $P$ is not irreducible as an $\mathbb{A}(k)$-module
if and only if
$k$ is singular, and $P/\operatorname{Rad}(k)$
is the unique simple quotient $L(triv,k)$ of $P$.
\end{prop}
\begin{proof}
Let $M\subset P$ be a proper $\mathbb{A}(k)$-submodule.
Because $1\not\in M$ we have, since $M$ is graded, $(1,M)_k=0$.
But this clearly implies $(P,M)_k=0$,
in other words: $M\subset \operatorname{Rad}(k)$.
Thus $\operatorname{Rad}(k)$ is the unique maximal
proper $\mathbb{A}(k)$-submodule of $P$.
\end{proof}

We are in the position to apply the technique of the Jantzen filtration
(see \cite{J}, Chapter 5), so let us discuss this briefly.
Given a real parameter $k_0\in K$, consider the real line $L$ in $K$
through $k_0$ and $0$.
Note that $k=0$ is
a regular parameter, hence a generic point on $L$ will be a regular parameter as well.
We parameterize this line by $k(t):=(1+t)k_0$ and consider $t$ as a real indeterminate.
Let us denote $R=\mathbb{R}[t]$.
For any complex vector space $B$ we denote by $B_R:=R\otimes_\mathbb{R} B$ the free
$R_c=\mathbb{C}\otimes_{\mathbb{R}} R=\mathbb{C}[t]$-module that arises
from $B$ by extension of scalars. Then $\mathbb{A}_R$ is
a free associative $R_c$-algebra.
For $r\in R_c$ we define $r^*(t):=\overline{r(t)}$ (recall that $\overline{t}=t$).
Thus $*$ is a $\mathbb{C}$-anti-linear involution on $R_c$.
We have that $P_R$ is a module for $\mathbb{A}_R$.
The anti-linear isomorphism $*:P\to S$ extends naturally to an anti-linear
isomorphism $*:P_R\to S_R$, where anti-linear means that $(rp)^*=r^*p^*$.
We define a contravariant hermitian form $(\cdot,\cdot)_R$ with values in $R_c$ on
$P_R$ by $(p,q)_R:=(p^*(T(k(t))q)(0)$. The form is linear in the second factor and
anti-linear in the first factor. It is hermitian in the sense that
$(p,q)_R=(q,p)_R^*$, and contravariant in the sense that
$(xp,q)_R=(p,T_{x^*}(k(t))q)_R$, $(T_{x^*}(k(t))p,q)_R=(p,xq)_R$ and
$(p^g,q^g)_R=(p,q)_R$.
Let $m_0$ denote the $*$-invariant maximal ideal $m_0=tR_c$ of $R_c$.
We introduce a sequence of $R_c$-linear subspaces $M^i_R\subset P_R$ defined by:
\begin{equation}
M^i_R:=\{p\in P_R \mid (P_R,p)_R\subset m_0^i\},
\end{equation}
\begin{lem}
The $M^i_R$ form a decreasing sequence of $\mathbb{A}_R$-submodules in $P_R$.
\end{lem}
\begin{proof}
The sequence $M^i_R$ is clearly a decreasing sequence of $R_c$-linear subspaces.
The fact that they are submodules follows from the contravariance of $(\cdot,\cdot)_R$.
\end{proof}
Let us denote by $\psi$ the specialization functor $\psi(L):=L/tL$ (where $L$ is an
$R_c$-module) at $t=0$. We thus obtain
a $\mathbb{C}$-algebra homomorphism $\psi:\mathbb{A}_R\to\mathbb{A}(k_0)$. This is
compatible with the module structures of $P_R$ and $P$,
in the sense that for $a\in\mathbb{A}_R$
and $p\in P_R$, we have $\psi(ap)=\psi(a)\psi(p)$. Also,
$\psi((p,q)_R)=(\psi(p),\psi(q))_{k_0}$. We put
\begin{equation}
M^i:=\psi(M^i_R).
\end{equation}
For a graded subspace $M$ of $P$ we introduce the notation
\begin{equation}
\operatorname{Ch}(M):=\sum_{n\in\mathbb{Z}_+}\operatorname{dim}(M_n)X^n.
\end{equation}
\begin{prop}
\begin{enumerate}
\item[(i)] $P_R=M^0\supset M^1\supset M^2\dots$ is a sequence of $\mathbb{A}(k_0)$-submodules.
\item[(ii)] We have $\sum_{i>0}\operatorname{Ch}(M^i)=
\sum_{n\geq 0}\nu(D(n))X^n$, where $D(n)$
denotes the determinant of $(\cdot,\cdot)_{R}$ on $P_{R,n}$ (the $R_c$-module
of polynomials of homogeneous degree $n$ on $V$, with coefficients in $R_c$),
and where $\nu$ denotes the $m_0$-adic valuation on the polynomial ring $R_c$.
\item[(iii)] For $i\gg 0$, $M^i=0$.
\item[(iv)] For $p\in M^i$, let us denote by $\tilde{p}\in M^i_R$ an arbitrary element
such that $\psi(\tilde{p})=p$.
For all $p,q\in M^i$, the expression
$(p,q)^{(i)}:=\psi(t^{-i}(\tilde{p},\tilde{q})_{R})$ depends only on $p$ and $q$,
not on the chosen lifts $\tilde{p},\tilde{q}$. The form $(\cdot,\cdot)^{(i)}$
is hermitian and $\mathbb{A}(k_0)$-contravariant on $M^i$,
and its radical is $M^{i+1}$.
In particular, $M^1=\operatorname{Rad}(k_0)$.
\end{enumerate}
\end{prop}
\begin{proof}
(i) This follows directly from the above Lemma, by application of $\psi$ to
the sequence $M^i_R$.

(ii) Since the spaces $P_{R,n}$ are mutually orthogonal with respect to
$(\cdot,\cdot)_R$, we see that
$M^i_{R,n}=\{p\in P_{R,n} \mid (P_{R,n},p)_{R}\in m_0^i\}$.
We apply (an appropriately adapted version of) Lemma 5.1 of \cite{J} to the
$R_c$-module $P_{R,n}$, and obtain that for all $n$
\begin{equation}
\sum_{i>0}\operatorname{dim}M^i_n=\nu(D(n)),
\end{equation}
where $M^i_n:=\psi(M^i_{R,n})=\psi(M^i_R)_n$.
Note that $D(n)$ is determined up to a unit in $R_c$, so that
$\nu(D(n))$ is well defined. The result follows.

(iii) By (ii) we see that given $n\in\mathbb{N}$, there exists a $b>0$ such
that $M^i_m=0$ for all $m<n$ and $i> b$. Take
$n>-\operatorname{max}_{\tau\in \hat{G}}c_\tau(k_0)$ and suppose that
$M^i\not=0$ for some $i> b$. Let $l\in\mathbb{N}$ be minimal such that
$M^i_{l}\not=0$, and let $\sigma\in \hat{G}$ be such that $M_{l,\sigma}^i\not=0$.
Then $T_\xi(k_0)M_{l,\sigma}^i=0$ for each $\xi$, and thus $E(k_0)M_{l,\sigma}^i=0$.
On the other hand, the eigenvalue of $E(k_0)$ on $M_{l,\sigma}^i$ equals
$l+c_\tau(k_0)>0$, since $l\geq n>-c_\tau(k_0)$. This is a contradiction.

(iv) This follows by \cite{J}, Bemerkung after Lemma 5.1. Notice that
the expression $(p,q)^{(i)}$ is independent of the lifts $\tilde{p}$ and
$\tilde{q}$ since, if $a\in M^i_R\cap tP$, we have
\begin{equation}
(a,M^i_R)_R\subset t(P,M^i_R)_R\subset m_0^{i+1}.
\end{equation}
When $i=0$, we have $(\cdot,\cdot)^{(0)}=(\cdot,\cdot)_{k_0}$. Thus the
result $M^1=\operatorname{Rad}(k_0)$ is the special case $i=0$.
\end{proof}
We end the section with a hint for the interpretation of the
multiplicities $\delta_\tau(k)=[P=M(triv,k),L(\tau,k)]$.
Suppose that $G$ has a Coxeter-like presentation (in the sense of \cite{BMR}, Appendix 2)
such that its diagram also provides a presentation of the fundamental group (the braid
group) of the regular orbit space of $G$. Suppose that the cyclotomic Hecke algebra
$H(G,u)$ corresponding to the diagram of $G$ can be generated by $|G|$ elements over the
ring $\mathbb{Z}[u,u^{-1}]$. The inhomogeneous relations for the simple generators
of $H(G,q)$ are of the form
\begin{equation}
(T_s-u_{C,0})(T_s-u_{C,1})\dots(T_s-u_{C,e_C-1})=0
\end{equation}
where $s$ is a reflection in a certain cyclic group $G_H$ with $H\in C\in\mathcal{C}$,
with determinant $\operatorname{det}(s)=\exp(2 \pi i/e_C)$. Let us now
view the parameter value $u_{C,j}$ as a function of the parameters $k_{C,j}$
as follows
\begin{equation}
u_{C,j}\to q_{C,j}:=\exp(2 \pi i(j-e_Ck_{C,j})/e_C).
\end{equation}
We extend the ring of definition of $H(G,q)$ to the ring $R$ of entire functions
in the parameters $k_{C,j}$ via this substitution. The resulting algebra is
denoted by $H(G)_R$.
It is known that the Hecke algebra $H(G)_K$ over the quotient field $K$ of $R$
is split semi-simple
(see \cite{O3}, Corollary 6.6), so that we can uniquely parameterize the irreducible
representations $\pi_\tau$ of $H(G)_K$ by the irreducible representations $\tau$ of $G$.
Let $U$ denote the principal indecomposable block of the trivial representation
in the $m$-adic completion $H(G)_m$ of $H(G)_R$, and let $K_m$ denote the quotient field of
the ring $R_m$ of formal power series in $k$ centered at $v$.
The results and method of \cite{DjO} seem to suggest that the multiplicities
$\delta_\tau^\prime(v)$ defined by
\begin{equation}
K_m\otimes U =\oplus_{\tau\in\hat{G}}\delta^\prime_\tau(v)(K_m\otimes (\pi_{\tau})_m)
\end{equation}
are related to the multiplicities $\delta_\tau(v)$ if $v_{C,j}<0$ for all $C,j$.

In the next section we will turn to the study of the infinite family of
imprimitive groups. The set $K^{sing}$ will be described in detail for this
class of complex reflection groups.
\section{The groups of type $G\left(  m,p,N\right)  $}
\subsection{Introduction}

In this section we study the particular case of the complex reflection group called
$G\left(  m,p,N\right)  $, which is a finite subgroup of $U\left(  N\right)
$. Because of its close relation to the symmetric group it is possible to
perform a detailed analysis of the Dunkl operators
(constructed for real reflection groups in \cite{D1}), the pairing, and the
analogues of the nonsymmetric Jack polynomials. In fact, the special case
$G\left(  2,1,N\right)  $ is exactly the hyperoctahedral group (type $B$), and
the results of one of the authors (Dunkl \cite{D4},\cite{D5}) on type-$B$
polynomials motivate the methods used in this section. Some of the notation
used in the first section is changed here to a mode better suited to deal with
monomials and permutations. The fundamental objects are polynomials in
$x=\left(  x_{1},x_{2},\ldots,x_{N}\right)  \in\mathbb{C}^{N}$ (considered as
coordinate functions); the group is realized as a subgroup of the matrix group
$U\left(  N\right)  $ acting on the row vector $x$. For a multi-index
(composition) $\alpha\in\mathbb{N}_{0}^{N}$ let $\left|  \alpha\right|
=\sum_{i=1}^{N}\alpha_{i}$ and let $x^{\alpha}$ denote the monomial
$\prod_{i=1}^{N}x_{i}^{\alpha_{i}}$. To a permutation $w\in S_{N}$ (the
symmetric group on $N$ letters) associate an $N\times N$ permutation matrix
with $1$'s at the $\left(  w\left(  j\right)  ,j\right)  $ entries. The action
on $x$ is given by $\left(  xw\right)  _{i}=x_{w\left(  i\right)  }$; the
action on polynomials is $wp\left(  x\right)  =p\left(  xw\right)  $. Thus the
action on monomials is $w\left(  x^{\alpha}\right)  =x^{w\alpha}$ where
$\left(  w\alpha\right)  _{i}=\alpha_{w^{-1}\left(  i\right)  }$ (consider
$\alpha$ as a column vector). The symmetric group contains the transpositions
$\left(  i,j\right)  $, for $i\neq j$, defined by
\[
x\left(  i,j\right)  =\left(  x_{1},\ldots,\overset{i}{x_{j}},\ldots
,\overset{j}{x_{i}},\ldots\right)  .
\]

For a fixed $m=2,3,\ldots$ let $\eta=e^{2\pi\mathrm{i}/m},$ an $m^{th}$ root
of unity, then the complex reflection group $W$ of type $G\left(
m,1,N\right)  $ consists of the $N\times N$ permutation matrices with the
nonzero entries being powers of $\eta$. The group is generated by the
transpositions $\left(  i,i+1\right)  $, $1\leq i\leq N-1$, and by the complex
reflection $\tau_{N}$ where $\tau_{i}$ is defined by
\[
x\tau_{i}=\left(  x_{1},\ldots,\overset{i}{\eta x_{i}},\ldots\right)  ,
\]
for $1\leq i\leq N$. The powers $\tau_{i}^{s}$ are also reflections. Thus
$\tau_{i}^{s}x^{\alpha}=\eta^{s\alpha_{i}}x^{\alpha}$ for $\alpha\in
\mathbb{N}_{0}^{N}$. The symmetric group $S_{N}$ is obviously a subgroup of
$W$. The group $W$ acts on polynomials by $wp\left(  x\right)  =p\left(
xw\right)  $ for $w\in W$. There are some obvious commutation relationships:
(where $x_{i}$ denotes the multiplication operator)
\begin{align*}
x_{i}\tau_{i}  & =\eta^{-1}\tau_{i}x_{i},\\
\left(  i,j\right)  \tau_{i}  & =\tau_{j}\left(  i,j\right)  ,
\end{align*}
and the elements $\tau_{i}^{-s}\left(  i,j\right)  \tau_{i}^{s}$ are ordinary
(period 2) reflections in $W$. In terms of root vectors $v\neq0$, such a
reflection is given by $x\sigma_{v}=x-2\left(  \left\langle v,x\right\rangle
/||v||^{2}\right)  v$ (where the hermitian inner product is $\left\langle
x,y\right\rangle =\sum_{j=1}^{N}\overline{x}_{j}y_{j}$ and the norm
$\vert$%
$\vert$%
$x||=\left\langle x,x\right\rangle ^{1/2}$).

\begin{definition}
\label{r0def}For $1\leq i\leq N$ let $e_{i}=(0,\ldots,0,\overset{i}{1},0,\ldots)$
denote the standard unit vector of $\mathbb{C}^N$.
For $i\neq j$ and $0\leq l\leq m-1$,
let $v_{ij}^{\left(  l\right)  }=e_{i}-\eta^{-l}e_{j}$, and $R_{0}=\left\{
v_{ij}^{\left(  l\right)  }:1\leq i<j\leq N,0\leq l\leq m-1\right\}  .$
\end{definition}

For $v=v_{ij}^{\left(  l\right)  }$ the reflection $\sigma_{v}$ equals
$\tau_{i}^{-l}\left(  i,j\right)  \tau_{i}^{l}$.
Similarly $x\tau_{i}^{l}=x-\left(
1-\eta^{l}\right)\allowbreak  \left\langle e_{i},x\right\rangle e_{i}$. To define the
Dunkl operators introduce $m$ parameters $\kappa_{s},0\leq
s\leq m-1$. In the notation of the Section 2, the class $C_{0}$ is
$\left\{  v^{\bot}:v\in R_{0}\right\}  $ and ${C}_{1}=\cup_{i=1}%
^{N}\left\{  x:x_{i}=0\right\}  $, further $e_{C_{0}}=2,\,k_{C_{0}}=\kappa
_{0}$ and $e_{C_{1}}=m,k_{C_{1},s}=\kappa_{s} $ for $1\leq s\leq m-1$. The
functional $\alpha_{H}=\left\langle v,\cdot\right\rangle $ where the
hyperplane $H=v^{\bot}$.

\begin{definition}
For $1\leq i\leq N$ let
\[
T_{i}=\frac{\partial}{\partial x_{i}}+\kappa_{0}\sum_{j\neq i}\sum_{s=0}%
^{m-1}\frac{1-\tau_{i}^{-s}\left(  i,j\right)  \tau_{i}^{s}}{x_{i}-\eta
^{s}x_{j}}+\sum_{t=1}^{m-1}\kappa_{t}\sum_{s=0}^{m-1}\frac{\eta^{-st}\tau
_{i}^{s}}{x_{i}},
\]
where the divisions are understood to follow the numerator operations.
\end{definition}

Note that $T_{i}$ is the same object as $T_{e_{i}}\left(  k\right)  $; the
nature of the group $W$ makes it desirable to use coordinates. The action can
be expressed in another useful way. Define associated projections on
polynomials by the linear extension of
\[
\pi_{j}\left(  s\right)  x^{\alpha}=\left\{
\begin{tabular}
[c]{l}%
$x^{\alpha}$ if $\alpha_{j}\equiv s\,\operatorname{mod}m$\\
$0$ else
\end{tabular}
\right.  ,
\]
for $1\leq j\leq N$ and $s\in\mathbb{N}_{0}$ . The projections can be
expressed as $\pi_{j}\left(  s\right)  =\varepsilon_{H_{j},s}=\frac{1}{m}%
\sum_{i=0}^{m-1}\eta^{-si}\tau_{j}^{i}$ (in the section 2 notation, where
$H_{i}=\left\{  x:x_{i}=0\right\}  $).

\begin{proposition}
\label{altdeft}For $\alpha\in\mathbb{N}_{0}^{N}$ and for $1\leq i\leq N$:
\begin{align*}
T_{i}x^{\alpha}  & =\frac{\partial}{\partial x_{i}}x^{\alpha}+m\kappa_{0}%
\sum_{j\neq i}\pi_{j}\left(  \alpha_{j}\right)  \frac{x^{\alpha}-\left(
i,j\right)  x^{\alpha}}{x_{i}-x_{j}}\\
& +\left\{
\begin{tabular}
[c]{l}%
$m\kappa_{s}x^{\alpha}/x_{i}$ if $\alpha_{i}\equiv s\,\operatorname{mod}m$ and
$1\leq s\leq m-1$\\
$0$ if $\alpha_{i}\equiv0\,\operatorname{mod}m$.
\end{tabular}
\right.  .
\end{align*}
\end{proposition}

\begin{proof}
The part involving (period 2) reflections ($\kappa_{0}$) is proven using the
following formulae: (stated for $i=1,j=2$, which suffices)
\begin{equation}
\frac{x_{1}^{\alpha_{1}}x_{2}^{\alpha_{2}}-\left(  1,2\right)  x_{1}%
^{\alpha_{1}}x_{2}^{\alpha_{2}}}{x_{1}-x_{2}}=\mathrm{sign}\left(  \alpha
_{1}-\alpha_{2}\right)  \sum_{t=\min\left(  \alpha_{1},\alpha_{2}\right)
}^{\max\left(  \alpha_{1},\alpha_{2}\right)  -1}x_{1}^{\alpha_{1}+\alpha
_{2}-t-1}x_{2}^{t},\label{quotij}%
\end{equation}
and for any $s$ with $0\leq s\leq m-1$, we have
\begin{align*}
\frac{x_{1}^{\alpha_{1}}x_{2}^{\alpha_{2}}-\tau_{1}^{-s}\left(  1,2\right)
\tau_{1}^{s}x_{1}^{\alpha_{1}}x_{2}^{\alpha_{2}}}{x_{1}-\eta^{s}x_{2}}  &
=\eta^{-\alpha_{2}s}\frac{x_{1}^{\alpha_{1}}\left(  \eta^{s}x_{2}\right)
^{\alpha_{2}}-x_{1}^{\alpha_{1}}\left(  \eta^{s}x_{2}\right)  ^{\alpha_{1}}%
}{x_{1}-\eta^{s}x_{2}}\\
& =\mathrm{sign}\left(  \alpha_{1}-\alpha_{2}\right)  \sum_{t=\min\left(
\alpha_{1},\alpha_{2}\right)  }^{\max\left(  \alpha_{1},\alpha_{2}\right)
-1}x_{1}^{\alpha_{1}+\alpha_{2}-t-1}x_{2}^{t}\eta^{s\left(  t-\alpha
_{2}\right)  };
\end{align*}
now the sum over $0\leq s\leq m-1$ in effect applies $m\pi_{2}\left(
\alpha_{2}\right)  $ to the sum in formula \ref{quotij}.
\end{proof}

\begin{remark}
\label{GmpN}The complex reflection group $G\left(  m,p,N\right)  $ (defined
when $p$ divides $m$) contains the reflections $\tau_{j}^{-s}\left(
i,j\right)  \tau_{j}^{s}$ (for $i<j$ and $0\leq s\leq m-1$) and $\tau_{i}%
^{sp}$ for $1\leq i\leq N,\,1\leq s\leq\frac{m}{p}-1$. The appropriate
modification of $\left\{  T_{i}\right\}  _{i=1}^{N}$ is to require
$\kappa_{sm/p}=0$ and $\kappa_{t+sm/p}=\kappa_{t}$ for $1\leq s\leq p-1$ and
$1\leq t\leq\frac{m}{p}-1$.
\end{remark}

\begin{proof}
Let $c_{s}=\sum_{j=1}^{m-1}\kappa_{j}\eta^{-sj}$ for $0\leq s\leq m-1$, then
$c_{s}$ is the coefficient of $\tau_{i}^{s}$ in the formula for $T_{i}$ (for
any $i$), and $1\leq s\leq m-1$. The inversion formula is $\kappa_{j}=\frac
{1}{m}\sum_{s=0}^{m-1}c_{s}\eta^{sj}$ for $j\geq1$, and $\sum_{s=0}^{m-1}%
c_{s}=0$. The condition $c_{s}=0$ unless $s\equiv0\operatorname{mod}p$ is
equivalent to the periodicity condition on the values of $\kappa_{j}$ stated above.
\end{proof}

We refer to the list of residues $\operatorname{mod}m$ of the index $\alpha$
as the parity type and say that $x^{\alpha}$ and $x^{\beta}$ have the same
parity type if $\alpha_{i}\equiv\beta_{i}\,\operatorname{mod}m$ for $1\leq
i\leq N$. If each monomial in a polynomial has the same parity type then we
say the polynomial has that type. Proposition \ref{altdeft} implies that if a
polynomial $p\left(  x\right)  $ has the same type as $x^{\alpha}$ for some
$\alpha\in\mathbb{N}_{0}^{N}$ then $T_{i}p\left(  x\right)  $ has the same
type as $x^{\alpha}x_{i}^{-1}$ (or $x^{\alpha}x_{i}^{m-1}$ if $\alpha_{i}=0$).
Thus the operators $\left\{  T_{i}x_{i}\right\}  _{i=1}^{N}$ preserve parity
type. We will define an inner-product structure on polynomials in which these
are self-adjoint. Just as in the symmetric group case they do not commute but
can be modified to form a commutative set. First we need to calculate the
commutant $\left[  x_{j},T_{i}\right]  =x_{j}T_{i}-T_{i}x_{j}$.

\begin{proposition}
For $i\neq j,\,\left[  x_{j},T_{i}\right]  =\kappa_{0}\sum_{s=0}^{m-1}\eta
^{s}\tau_{j}^{-s}\left(  i,j\right)  \tau_{j}^{s}$.
\end{proposition}

\begin{proof}
Use the definition of $T_{i}$ and consider $x_{j}T_{i}-T_{i}x_{j}$. The terms
involving differentiation and transpositions $\left(  i,t\right)  $ with
$t\neq i,j$ cancel, leaving only
\begin{align*}
x_{j}T_{i}-T_{i}x_{j}  & =\kappa_{0}\sum_{s=0}^{m-1}\left(  x_{j}\frac
{1-\tau_{i}^{-s}\left(  i,j\right)  \tau_{i}^{s}}{x_{i}-\eta^{s}x_{j}}%
-\frac{x_{j}-\tau_{i}^{-s}\left(  i,j\right)  \tau_{i}^{s}x_{j}}{x_{i}%
-\eta^{s}x_{j}}\right) \\
& =\kappa_{0}\sum_{s=0}^{m-1}\frac{x_{i}\eta^{-s}-x_{j}}{x_{i}-\eta^{s}x_{j}%
}\tau_{i}^{-s}\left(  i,j\right)  \tau_{i}^{s}=\kappa_{0}\sum_{s=0}^{m-1}%
\eta^{-s}\tau_{i}^{-s}\left(  i,j\right)  \tau_{i}^{s}.
\end{align*}
This completes the proof.
\end{proof}

\begin{corollary}
\label{tixitjxj}For $i\neq j$ the following hold:
\begin{gather*}
\left[  T_{i}x_{i},T_{j}x_{j}\right]  =\kappa_{0}\left(  T_{i}x_{i}-T_{j}%
x_{j}\right)  \sum_{s=0}^{m-1}\tau_{j}^{-s}\left(  i,j\right)  \tau_{j}^{s},\\
\left[  T_{i}x_{i},T_{j}x_{j}-\kappa_{0}\sum\nolimits_{s=0}^{m-1}\tau_{j}%
^{-s}\left(  i,j\right)  \tau_{j}^{s}\right]  =0.
\end{gather*}
\end{corollary}

\begin{proof}
Indeed
\[
T_{i}x_{i}T_{j}x_{j}-T_{j}x_{j}T_{i}x_{i}=T_{i}\left(  x_{i}T_{j}-T_{j}%
x_{i}\right)  x_{j}-T_{j}\left(  x_{j}T_{i}-T_{i}x_{j}\right)  x_{i},
\]
and by the Proposition
\begin{align*}
T_{i}\left(  x_{i}T_{j}-T_{j}x_{i}\right)  x_{j}  & =\kappa_{0}T_{i}\sum
_{s=0}^{m-1}\eta^{s}\tau_{i}^{-s}\left(  i,j\right)  \tau_{i}^{s}x_{j}\\
& =\kappa_{0}T_{i}\sum_{s=0}^{m-1}\eta^{s}\tau_{i}^{-s}x_{i}\left(
i,j\right)  \tau_{i}^{s}\\
& =\kappa_{0}T_{i}x_{i}\sum_{s=0}^{m-1}\tau_{i}^{-s}\left(  i,j\right)
\tau_{i}^{s}%
\end{align*}
and similarly $T_{j}\left(  x_{j}T_{i}-T_{i}x_{j}\right)  x_{i}=\allowbreak
\kappa_{0}T_{j}x_{j}\sum_{s=0}^{m-1}\tau_{j}^{-s}\left(  i,j\right)  \tau
_{j}^{s}.$ Finally observe that $\sum_{s=0}^{m-1}\tau_{i}^{-s}\left(  i,j\right)
\tau_{i}^{s}=\sum_{s=0}^{m-1}\tau_{i}^{-s}\tau_{j}^{s}\left(  i,j\right)
\allowbreak=\allowbreak\sum_{s=0}^{m-1}\tau_{i}^{s}\tau_{j}^{-s}\left(  i,j\right)
\allowbreak =\sum_{s=0}^{m-1}\tau_{j}^{-s}\left(  i,j\right)  \tau_{j}^{s}$ (by changing
the index of summation from $s$ to $-s$ and using the relation $\tau_{i}%
\tau_{j}=\tau_{j}\tau_{i}$).

Furthermore the commutant satisfies $\left[  T_{i}x_{i},\tau_{j}^{-s}\left(  i,j\right)
\tau_{j}^{s}\right]  =T_{i}x_{i}\tau_{j}^{-s}\left(  i,j\right)  \tau_{j}%
^{s}-\tau_{j}^{-s}\left(  i,j\right)  \tau_{j}^{s}T_{i}x_{i}=T_{i}x_{i}%
\tau_{j}^{-s}\left(  i,j\right)  \tau_{j}^{s}-T_{j}x_{j}\tau_{j}^{-s}\left(
i,j\right)  \tau_{j}^{s}$ for each $s $, since $\tau_{j}$ commutes with
$T_{j}x_{j}$.
\end{proof}

\begin{definition}
For $1\leq i\leq N$ define the operators
\[
\mathcal{U}_{i}=T_{i}x_{i}-\kappa_{0}\sum_{j<i}\sum_{s=0}^{m-1}\tau_{i}%
^{-s}\left(  i,j\right)  \tau_{i}^{s}.
\]
\end{definition}

\begin{theorem}
$\mathcal{U}_{i}\mathcal{U}_{j}=\mathcal{U}_{j}\mathcal{U}_{i}.$
\end{theorem}

\begin{proof}
Let $\lambda_{rt}=\sum_{s=0}^{m-1}\tau_{r}^{-s}\left(  r,t\right)  \tau
_{r}^{s}=\sum_{s=0}^{m-1}\tau_{t}^{-s}\left(  r,t\right)  \tau_{t}^{s}$.
Suppose that $i<j$, then we have $\left[  \mathcal{U}_{i},\mathcal{U}_{j}\right]  =\left[
T_{i}x_{i},T_{j}x_{j}-\kappa_{0}\lambda_{ij}\right]  -\allowbreak\kappa
_{0}\sum_{r<j,r\neq i}\left[  T_{i}x_{i},\lambda_{rj}\right]  -\kappa_{0}%
\sum_{r<i}\left[  \lambda_{ri},T_{j}x_{j}\right]  +\allowbreak\kappa_{0}%
^{2}\sum_{r<i}\left(  \left[  \lambda_{ri},\lambda_{rj}\right]  +\left[
\lambda_{ri},\lambda_{ij}\right]  \right)  .$ The commutators $\left[
\lambda_{ri},\lambda_{tj}\right]  $ with $t\neq r,i$ vanish. In the previous
expansion all terms but the last are zero (by the Corollary). For $r<i<j$
(this computation is in the group algebra $\mathbb{C}W$)
\begin{align*}
\lambda_{ri}\lambda_{rj}  & =\sum_{s=0}^{m-1}\tau_{r}^{-s}\left(  r,i\right)
\tau_{r}^{s}\sum_{t=0}^{m-1}\tau_{j}^{-t}\left(  r,j\right)  \tau_{j}^{t}%
=\sum_{s,t=0}^{m-1}\tau_{r}^{-s}\tau_{i}^{s+t}\tau_{j}^{-t}\left(  r,i\right)
\left(  r,j\right) \\
& =\sum_{s,t=0}^{m-1}\tau_{r}^{-s}\tau_{i}^{s+t}\tau_{j}^{-t}\left(
i,j\right)  \left(  r,i\right)  =\sum_{t=0}^{m-1}\tau_{j}^{-t}\left(
i,j\right)  \tau_{j}^{t}\sum_{s=0}^{m-1}\tau_{r}^{-s}\left(  r,i\right)
\tau_{r}^{s}\\
& =\lambda_{ij}\lambda_{ri},
\end{align*}
and similarly $\lambda_{ri}\lambda_{ij}=\lambda_{rj}\lambda_{ri}$. Thus
$\left[  \lambda_{ri},\lambda_{rj}\right]  +\left[  \lambda_{ri},\lambda
_{ij}\right]  =0$.
\end{proof}

The group algebra elements $\sum_{r<i}\lambda_{ri}$ (for $1<i\leq N$) are the
analogues of the Jucys-Murphy elements for the symmetric group. The inner
product we will use is the pairing $\left(  p,q\right)  _{k}=p^{\ast}\left(
T\right)  q\left(  x\right)  |_{x=0}$, which means that the operator $p^{\ast
}\left(  T\right)  ,$ obtained from $p^{\ast}\left(  x\right)  =\overline
{p\left(  \overline{x}\right)  }$ (each coefficient of $p$ is conjugated) by
replacing each $x_{i}$ by $T_{i}$, is applied to the polynomial $q\left(
x\right)  $ and the resulting polynomial is evaluated at $x=0$. Obviously if
$p,q$ are homogeneous of the same degree then $p^{\ast}\left(  T\right)
q\left(  x\right)  $ is a constant (degree $0)$; and if $p,q$ are homogeneous
of different degrees then $\left(  p,q\right)  _{k}=0$. Of course this is the
coordinatized form of the general definition given in \ref{eq:pair}. The
known results for $S_{N}$ are used to analyze the action of $T_{i}$ by
factoring polynomials into a ``parity part'' and a part invariant under the
complex reflections $\tau_{i}.$ Introduce a new variable
\[
y=\left(  y_{1},\ldots,y_{N}\right)  =\left(  x_{1}^{m},\ldots,x_{N}%
^{m}\right)  .
\]
Say that a composition $\alpha\in\mathbb{N}_{0}^{N}$ is a \textit{parity type}
if $0\leq\alpha_{i}<m$ for each $i$. If all the terms of a polynomial
$p\left(  x\right)  $ have the same parity type then $p$ can be expressed in
the form $p\left(  x\right)  =x^{\alpha}g\left(  y\right)  $ for some parity
type $\alpha.$ The type-$A$ Dunkl operators (see
Section 8.3 in the monograph \cite{DX} by Dunkl and Xu) are defined by
\[
\mathcal{D}_{i}g\left(  y\right)  =\frac{\partial}{\partial y_{i}}g\left(
y\right)  +\kappa_{0}\sum_{j\neq i}\frac{g\left(  y\right)  -\left(
i,j\right)  g\left(  y\right)  }{y_{i}-y_{j}},
\]
for $1\leq i\leq N$ (using the same notation $\left(  i,j\right)  $ for the
transpositions acting on $y$ as on $x$). Often we write $\frac{1-\left(
i,j\right)  }{y_{i}-y_{j}}g\left(  y\right)  $ for the inner term.

\begin{proposition}
\label{tixg}Suppose $g\left(  y\right)  $ is a polynomial and $\alpha$ is a
parity type then

\begin{enumerate}
\item  if $\alpha_{i}>0$ then\newline $T_{i}x^{\alpha}g\left(  y\right)
=mx^{\alpha}x_{i}^{-1}\left(  \dfrac{\alpha_{i}}{m}+\kappa_{\alpha_{i}%
}+\mathcal{D}_{i}y_{i}-1-\kappa_{0}\sum\limits_{j\neq i,\alpha_{j}\geq
\alpha_{i}}\left(  i,j\right)  \right)  g\left(  y\right)  ;$

\item  if $\alpha_{i}=0$ then $T_{i}x^{\alpha}g\left(  y\right)  =mx^{\alpha
}x_{i}^{m-1}\mathcal{D}_{i}g\left(  y\right)  .$
\end{enumerate}
\end{proposition}

\begin{proof}
Suppose $\alpha_{i}>0$ then $T_{i}x^{\alpha}g\left(  y\right)  =\left(
\alpha_{i}+m\kappa_{\alpha_{i}}\right)  x^{\alpha}x_{i}^{-1}\allowbreak +
mx^{\alpha}x_{i}^{m-1}\frac{\partial}{\partial y_{i}}g\left(  y\right)\allowbreak
+m\kappa_{0}\sum\limits_{j\neq i}E_{j}$, where (for $j\neq i)$
\begin{align*}
E_{j}  & =\pi_{j}\left(  \alpha_{j}\right)  \dfrac{x^{\alpha}g\left(
y\right)  -\left(  i,j\right)  x^{\alpha}g\left(  y\right)  }{x_{i}-x_{j}}\\
& =\pi_{j}\left(  \alpha_{j}\right)  \left(  x^{\alpha}\dfrac{y_{i}-y_{j}%
}{x_{i}-x_{j}}\dfrac{1-\left(  i,j\right)  }{y_{i}-y_{j}}g\left(  y\right)
+\dfrac{x^{\alpha}-\left(  i,j\right)  x^{\alpha}}{x_{i}-x_{j}}\left(
i,j\right)  g\left(  y\right)  \right)  .
\end{align*}
For the first term note $\pi_{j}\left(  \alpha_{j}\right)  x^{\alpha}%
\dfrac{y_{i}-y_{j}}{x_{i}-x_{j}}=\pi_{j}\left(  \alpha_{j}\right)  x^{\alpha
}\sum\limits_{s=0}^{m-1}x_{i}^{m-1-s}x_{j}^{s}=x_{i}^{m-1}x^{\alpha}$. By
formula \ref{quotij} (with $1,2$ replaced by $i,j$)
\[
\dfrac{x^{\alpha}-\left(  i,j\right)  x^{\alpha}}{x_{i}-x_{j}}=\mathrm{sign}%
\left(  \alpha_{i}-\alpha_{j}\right)  \sum_{s=\min\left(  \alpha_{i}%
,\alpha_{j}\right)  }^{\max\left(  \alpha_{i},\alpha_{j}\right)  -1}x^{\alpha
}x_{i}^{\alpha_{j}-s-1}x_{j}^{s-\alpha_{j}},
\]
(note that the power of $x_{j}$ is $s$) and it follows that $\pi_{j}\left(
\alpha_{j}\right)  \dfrac{x^{\alpha}-\left(  i,j\right)  x^{\alpha}}%
{x_{i}-x_{j}}=x^{\alpha}x_{i}^{-1}$ if $\alpha_{i}>\alpha_{j}$ and $=0$ if
$\alpha_{i}\leq\alpha_{j}$; if $\alpha_{i}>\alpha_{j}$ then the summation
extends over $0\leq\alpha_{j}\leq s\leq\alpha_{i}-1\leq m-2$ and the
projection $\pi_{j}\left(  \alpha_{j}\right)  $ picks out the term with
$s=\alpha_{j}$, and when $\alpha_{i}\leq\alpha_{j}$ the case $s=\alpha_{j}$
cannot occur. Also $\mathcal{D}_{i}y_{i}g\left(  y\right)  =g\left(  y\right)
+y_{i}\frac{\partial}{\partial y_{i}}g\left(  y\right)  +\kappa_{0}%
\sum\limits_{j\neq i}\dfrac{y_{i}g\left(  y\right)  -y_{j}\left(  i,j\right)
g\left(  y\right)  }{y_{i}-y_{j}}$ and the inner term equals $y_{i}%
\dfrac{1-\left(  i,j\right)  }{y_{i}-y_{j}}g\left(  y\right)  +\left(
i,j\right)  g\left(  y\right)  $. Consequently, for each $j\neq i$ we obtain
the equality $x^{\alpha}%
x_{i}^{m-1}\dfrac{1-\left(  i,j\right)  }{y_{i}-y_{j}}g\left(  y\right)
=x^{\alpha}x_{i}^{-1}y_{i}\dfrac{1-\left(  i,j\right)  }{y_{i}-y_{j}}g\left(
y\right)  =\allowbreak x^{\alpha}x_{i}^{-1}\left(  \dfrac{1-\left(
i,j\right)  }{y_{i}-y_{j}}y_{i}g\left(  y\right)  -\left(  i,j\right)
g\left(  y\right)  \right)  .$ The transpositions $\left(  i,j\right)  $ for
$\alpha_{i}>\alpha_{j}$ are cancelled out.

If $\alpha_{i}=0$ then $T_{i}x^{\alpha}g\left(  y\right)  =mx^{\alpha}%
x_{i}^{m-1}\frac{\partial}{\partial y_{i}}g\left(  y\right)  +m\kappa_{0}%
\sum\limits_{j\neq i}E_{j}$ with the same $E_{j}$ as before. The case
$\alpha_{i}>\alpha_{j}$ can not occur so that $E_{j}=x_{i}^{m-1}x^{\alpha
}\dfrac{1-\left(  i,j\right)  }{y_{i}-y_{j}}g\left(  y\right)  $. This
completes the proof.
\end{proof}

\subsection{Simultaneous Eigenfunctions}

We use the nonsymmetric Jack polynomials from the type-$A$ machinery to
produce a complete set of simultaneous eigenfunction for the commuting operators
$\left\{  \mathcal{U}_{i}:1\leq i\leq N\right\}  $. The nicest case is for
polynomials of ``\textit{standard parity type}'', that is, of the form
$x^{\alpha}g\left(  y\right)  $ such that $i<j$ implies $\alpha_{i}\geq
\alpha_{j}$. This means that in the list $\left(  \alpha_{1},\alpha_{2}%
,\ldots,\alpha_{N}\right)  $ the values $\alpha_{i}=m-1$ appear first, then
the values $\alpha_{i}=m-2$ and so on until the end of the list composed of
$\alpha_{i}=0$ values (roughly $\left(  m-1,\ldots,m-1,m-2,\ldots
,m-2,\ldots,0\right)  $). Not every value need appear, of course. First we
restate Proposition \ref{tixg} for the commuting operators.

\begin{proposition}
\label{uxg}Suppose $g\left(  y\right)  $ is a polynomial and $\alpha$ is a
parity type then
\begin{enumerate}
\item  if $0\leq\alpha_{i}<m-1$ then $\mathcal{U}_{i}x^{\alpha
}g\left(  y\right)=  \newline =mx^{\alpha}\left(  \dfrac{\alpha_{i}+1}{m}+\kappa
_{\alpha_{i}+1}+\mathcal{D}_{i}y_{i}-1-\kappa_{0}\left(  \sum\limits_{\alpha
_{j}>\alpha_{i}}+\sum\limits_{j<i,\alpha_{j}=\alpha_{i}}\right)  \left(
i,j\right)  \right)  g\left(  y\right)  ;$

\item  if $\alpha_{i}=m-1$ then $\mathcal{U}_{i}x^{\alpha}g\left(  y\right)
=mx^{\alpha}\left(  \mathcal{D}_{i}y_{i}-\kappa_{0}\sum\limits_{j<i,\alpha
_{j}=m-1}\left(  i,j\right)  \right)  g\left(  y\right)  .$
\end{enumerate}
\end{proposition}

\begin{proof}
For any $\beta\in\mathbb{N}_{0}^{N}$ we have $\sum_{s=0}^{m-1}\tau_{i}%
^{-s}\left(  i,j\right)  \tau_{i}^{s}x^{\beta}\allowbreak =\sum_{s=0}^{m-1}%
\eta^{s\left(  \beta_{i}-\beta_{j}\right)  }\left(  i,j\right)  x^{\beta
}\allowbreak =m\left(  i,j\right)  x^{\beta}$ if $\beta_{j}\equiv\beta_{i}%
\operatorname{mod}m$ and else equals $0$. Thus the second part of
$\mathcal{U}_{i}$ contributes $-\kappa_{0}\sum\limits_{j<i,\alpha_{j}%
=\alpha_{i}}\left(  i,j\right)  g\left(  y\right)  $ (and $\left(  i,j\right)
x^{\alpha}=x^{\alpha}$ for such $j$). When $\alpha_{i}=m-1$ then write
$x_{i}x^{\alpha}$ as $x^{\alpha}x_{i}^{1-m}y_{i}$ and use part (2) of
Proposition \ref{tixg}.
\end{proof}

Following Definition 8.3.4 in \cite{DX}, the type-$A$ (commuting) operators
are defined by
\[
\mathcal{U}_{i}^{A}=\mathcal{D}_{i}y_{i}+\kappa_{0}-\kappa_{0}\sum
\limits_{j<i}\left(  i,j\right)  .
\]
The nonsymmetric Jack polynomials are defined by means of a partial order on compositions.

\begin{definition}
\label{defdom}For $\mu,\nu\in\mathbb{N}_{0}^{N}$ the relation $\mu\succ\nu$
means $\sum_{i=1}^{j}\mu_{i}\geq\sum_{i=1}^{j}\nu_{i}\ $for each $j$ and
$\mu\neq\nu$ (\textit{dominance} order); $\mu^{+}\ $is defined to be the
(unique) partition $w\mu$ (for some $w\in S_{N}$; that is,$1\leq$ $i<j\leq N$
implies $\left(  \mu^{+}\right)  _{i}\geq\left(  \mu^{+}\right)  _{j}$ ); and
$\mu\vartriangleright\nu$ means $\left|  \mu\right|  =\left|  v\right|  $ and
$\mu^{+}\succ\nu^{+}$ or $\mu^{+}=\nu^{+}$and $\mu\succ\nu$.
\end{definition}

Define a convenient basis $\left\{  p_{\mu}:\mu\in\mathbb{N}_{0}^{N}\right\}
$ for homogeneous polynomials (see \cite{D3}) by the generating function:
\[
\sum_{\mu\in\mathbb{N}_{0}^{N}}p_{\mu}\left(  y\right)  z^{\mu}=\prod
_{i=1}^{N}\left\{  \left(  1-y_{i}z_{i}\right)  ^{-1}\prod_{j=1}^{N}\left(
1-y_{j}z_{i}\right)  ^{-\kappa_{0}}\right\}  ;
\]
with the useful property that $\mu_{i}=0$ implies $\mathcal{D}_{i}p_{\mu}=0$.
For each $\mu\in\mathbb{N}_{0}^{N}$ there is a unique simultaneous
eigenfunction of $\left\{  \mathcal{U}_{i}^{A}\right\}  _{i=1}^{N}$ of the
form $\zeta_{\mu}=p_{\mu}+\sum_{\nu\vartriangleright\mu}B\left(  v,\mu\right)
p_{\nu}$ (where $B\left(  v,\mu\right)  \in\mathbb{Q}\left(  \kappa
_{0}\right)  $ and does not depend on $N$ provided that $N\geq M$ and $\mu
_{i}=0$ for all $i>M$). See Theorem 8.4.13 in \cite{DX}.

Then $\zeta_{\mu}\left(  y\right)  $ satisfies $\mathcal{U}_{i}^{A}\zeta_{\mu
}\left(  y\right)  =\xi_{i}\left(  \mu\right)  \zeta_{\mu}\left(  y\right)  $
for $\mu\in\mathbb{N}_{0}^{N}$ and $1\leq i\leq N$; where the eigenvalues are
given by
\[
\xi_{i}\left(  \mu\right)  =\kappa_{0}\left(  N-\#\left\{  j:\mu_{j}>\mu
_{i}\right\}  -\#\left\{  j:j<i\,\&\,\mu_{j}=\mu_{i}\right\}  \right)
+\mu_{i}+1.
\]
Suppose $\beta\in\mathbb{N}_{0}^{N}$ then there is a unique parity type
$\alpha$ and a composition $\gamma\in\mathbb{N}_{0}^{N}$ so that $\beta
=\alpha+m\gamma$ (as vectors; that is, $\beta_{i}=\alpha_{i}+m\gamma_{i} $ and
$\gamma_{i}=\left\lfloor \beta_{i}/m\right\rfloor $ for each $i$). We will
construct a simultaneous eigenfunction for each $\beta$. When $\alpha$ is a
standard parity type the nonsymmetric Jack polynomials work directly.

\begin{proposition}
\label{uxg0}Suppose $\alpha$ is a standard parity type, $g$ is any polynomial
in $y$, and $1\leq i\leq N$, then

\begin{enumerate}
\item  if $0\leq\alpha_{i}<m-1$ then $\mathcal{U}_{i}x^{\alpha}g\left(
y\right)  =mx^{\alpha}\left(  \frac{\alpha_{i}+1}{m}+\kappa_{\alpha_{i}%
+1}-\kappa_{0}-1+\mathcal{U}_{i}^{A}\right)  g\left(  y\right)  ;$

\item  if $\alpha_{i}=m-1$ then $\mathcal{U}_{i}x^{\alpha}g\left(  y\right)
=mx^{\alpha}\left(  \mathcal{U}_{i}^{A}-\kappa_{0}\right)  g\left(  y\right)
. $
\end{enumerate}
\end{proposition}

\begin{proof}
By definition of standard parity type, for any $i$ the set $\left\{
j:\alpha_{j}>\alpha_{i}\right\}  \allowbreak\cup\left\{  j:j<i\,\&\,\alpha
_{j}=\alpha_{i}\right\}  $ is exactly the set $\left\{  j:1\leq j<i\right\}  .$
\end{proof}

\begin{corollary}
Suppose $\alpha$ is a standard parity type, and $\gamma\in\mathbb{N}_{0}^{N}$ then

\begin{enumerate}
\item $\mathcal{U}_{i}x^{\alpha}\zeta_{\gamma}\left(  y\right)  =m\left(
\frac{\alpha_{i}+1}{m}-1+\kappa_{\alpha_{i}+1}+\xi_{i}\left(  \gamma\right)
-\kappa_{0}\right)  x^{\alpha}\zeta_{\gamma}\left(  y\right)  $ when
$0\leq\alpha_{i}<m-1;$

\item $\mathcal{U}_{i}x^{\alpha}\zeta_{\gamma}\left(  y\right)  =m\left(
\xi_{i}\left(  \gamma\right)  -\kappa_{0}\right)  x^{\alpha}\zeta_{\gamma
}\left(  y\right)  $ when $\alpha_{i}=m-1.$
\end{enumerate}
\end{corollary}
Suppose $\alpha^{\prime}$ is any parity type, but not standard. The idea is to
use a permutation which maps $\alpha^{\prime}$ to a standard parity type in
such a way that the original order of coordinates with the same value of
$\alpha_{i}^{\prime}$ is preserved. For technical reasons it is easier to do
this backwards. Suppose that $\alpha$ is a standard parity type and suppose
$w\in S_{N}$ (a permutation) and has the property that $1\leq i<j\leq N$ and
$\alpha_{i}=\alpha_{j}$ implies $w\left(  i\right)  <w\left(  j\right)  $. We
will show that $wx^{\alpha}\zeta_{\beta}\left(  y\right)  $ is an eigenvector
of each $\mathcal{U}_{i}$. We use the transformation properties (holding for
any $w\in S_{N}$) : $w\mathcal{D}_{i}y_{i}=\mathcal{D}_{w\left(  i\right)
}y_{w\left(  i\right)  }w$ and $\left(  r,s\right)  w=w\left(  w^{-1}\left(
r\right)  ,w^{-1}\left(  s\right)  \right)  $ for $r\neq s$ (the action of
permutations is as follows: $w$ is an $N\times N $ permutation matrix with 1's
at the $\left(  w\left(  i\right)  ,i\right)  $ entries, $x$ is a row vector,
compositions are column vectors; so that $w\left(  x^{\alpha}\right)
=x^{w\alpha}$ where $\left(  w\alpha\right)  _{w\left(  i\right)  }=\alpha
_{i}$ for any $i$). In the following $wx^{\alpha}\zeta_{\gamma}\left(
y\right)  $ is the polynomial $x^{w\alpha}\zeta_{\gamma}\left(  yw\right)  $
(with parity type $w\alpha$).

\begin{proposition}
\label{nonstd}Suppose $\alpha$ is a standard parity type, $\gamma\in
\mathbb{N}_{0}^{N}$, and $w\in S_{N}$ has the property that $w\left(
i\right)  <w\left(  j\right)  $ whenever $1\leq i<j\leq N$ and $\alpha
_{i}=\alpha_{j}$, then

\begin{enumerate}
\item $\mathcal{U}_{w\left(  i\right)  }wx^{\alpha}\zeta_{\gamma}\left(
y\right)  =m\left(  \frac{\alpha_{i}+1}{m}-1+\kappa_{\alpha_{i}+1}+\xi
_{i}\left(  \gamma\right)  -\kappa_{0}\right)  wx^{\alpha}\zeta_{\gamma
}\left(  y\right)  $ when $0\leq\alpha_{i}<m-1;$

\item $\mathcal{U}_{w\left(  i\right)  }wx^{\alpha}\zeta_{\gamma}\left(
y\right)  =m\left(  \xi_{i}\left(  \gamma\right)  -\kappa_{0}\right)
wx^{\alpha}\zeta_{\gamma}\left(  y\right)  $ when $\alpha_{i}=m-1.$
\end{enumerate}
\end{proposition}

\begin{proof}
Suppose $\alpha_{i}<m-1$ then part (1) of Proposition \ref{uxg} applies (and
note that $\left(  w\alpha\right)  _{w\left(  i\right)  }=\alpha_{i}$) so
that
\begin{align*}
& \mathcal{U}_{w\left(  i\right)  }wx^{\alpha}\zeta_{\gamma}\left(  y\right)
\\
& =mx^{w\alpha}\left(  \dfrac{\alpha_{i}+1}{m}+\kappa_{\alpha_{i}%
+1}+\mathcal{D}_{w\left(  i\right)  }y_{w\left(  i\right)  }-1-\kappa_{0}%
\sum_{j\in E}\left(  w\left(  i\right)  ,j\right)  \right)  w\zeta_{\gamma}\\
& =mx^{w\alpha}w\left(  \dfrac{\alpha_{i}+1}{m}+\kappa_{\alpha_{i}%
+1}+\mathcal{D}_{i}y_{i}-1-\kappa_{0}\sum_{j\in E}\left(  i,w^{-1}\left(
j\right)  \right)  \right)  \zeta_{\gamma};
\end{align*}
where the set $E=\left\{  j:\left(  w\alpha\right)  _{j}>\left(
w\alpha\right)  _{w(i)}\right\}  \cup\left\{  j:j<w\left(  i\right)
\,\&\,\left(  w\alpha\right)  _{j}=\left(  w\alpha\right)  _{w(i)}\right\}  $.
But $\sum\limits_{j\in E}\left(  i,w^{-1}\left(  j\right)  \right)
=\allowbreak\sum\limits_{w\left(  s\right)  \in E}\left(  i,s\right)
=\sum\limits_{s\in w^{-1}E}\left(  i,s\right)  $, and $w^{-1}E=\left\{
s:\alpha_{s}>\alpha_{i}\right\}  \cup\left\{  s:w\left(  s\right)  <w\left(
i\right)  \,\&\,\alpha_{s}=\alpha_{i}\right\}  .$ By the condition on $w$, the
set $w^{-1}E$ is equal to $\left\{  s:1\leq s<i\right\}  $. Indeed, consider any $j<i$;
$\alpha_{j}<\alpha_{i}$ is impossible by definition of standard parity type so
$\alpha_{j}\geq\alpha_{i}$; furthermore if $\alpha_{j}=\alpha_{i}$ for some
$j$ then $j<i $ if and only if $w\left(  j\right)  <w\left(  i\right)  .$\ The
proof of part (1) is now completed similarly to the proof of Proposition
\ref{uxg0}.

The proof of part (2) is an obvious modification of that for part (1).
\end{proof}

\subsection{Evaluation of the Pairing}

We recall some facts from Proposition \ref{prop:easy} and Theorem \ref{thm:herm}.

\begin{theorem}
\label{pair}The pairing $\left(  \cdot,\cdot\right)  _{k}$ has the following
properties (with $p,q\in\mathcal{P}_{n}$):

\begin{enumerate}
\item $\left(  q,p\right)  _{k}=\overline{\left(  p,q\right)  _{k}}$ ;

\item $\left(  T_{i}x_{i}p,q\right)  _{k}=\left(  p,T_{i}x_{i}q\right)  _{k}$
\ for $1\leq i\leq N$;

\item $\left(  wp,wq\right)  _{k}=\left(  p,q\right)  _{k}$ \ for any $w\in W.$
\end{enumerate}
\end{theorem}

\begin{proof}
From the definition it is clear that $\left(  p,T_{i}x_{i}q\right)
_{k}=\left(  x_{i}p,x_{i}q\right)  _{k}$, and by part (1) we have $\left(
T_{i}x_{i}p,q\right)  _{k}=\overline{\left(  q,T_{i}x_{i}p\right)  _{k}%
}=\overline{\left(  x_{i}q,x_{i}p\right)  _{k}}=\left(  x_{i}p,x_{i}q\right)
_{k}. $

For the real case $(m=2$) there is an associated inner product with respect to
the measure $\prod_{i=1}^{N}\left|  x_{i}\right|  ^{2\kappa_{1}}\prod
_{i<j}\left|  x_{i}^{2}-x_{j}^{2}\right|  ^{2\kappa_{0}}\exp\left(
-\frac{\left|  x\right|  ^{2}}{2}\right)  dx$ on $\mathbb{R}^{N}\,$(see
\cite{D2} or Theorem 5.2.7 and Section 9.6.3 in \cite{DX}). It is an
interesting question whether there is a similar situation for the complex group.
\end{proof}

By the properties of the pairing in Theorem \ref{pair} there are obvious
orthogonality relations.

\begin{lemma}\label{lem:sadj}
The operators $\mathcal{U}_{i}$ are self-adjoint for $\left(  \cdot
,\cdot\right)  _{k}$.
\end{lemma}

\begin{proof}
The reflections $\sigma=\tau_{i}^{-s}\left(  i,j\right)  \tau_{i}^{s}$ are
self-adjoint (by part (2) of Theorem \ref{pair} because $\sigma^{2}=1$. By
part (3) of this theorem $T_{i}x_{i}$ is self-adjoint, thus, $\mathcal{U}_{i}
$ is also.
\end{proof}

\begin{proposition}
Suppose $\alpha,\beta$ are parity types and $\mu,\nu\in\mathbb{N}_{0}^{N}$ then:

\begin{enumerate}
\item  if $g_{1}\left(  y\right)  ,g_{2}\left(  y\right)  $ are polynomials
and $\alpha\neq\beta$ then $\left(  x^{\alpha}g_{1}\left(  y\right)
,x^{\beta}g_{2}\left(  y\right)  \right)  _{k}=0$;

\item  if $\mu\neq\nu$ and $\alpha$ is a standard parity type then $\left(
x^{\alpha}\zeta_{\mu}\left(  y\right)  ,x^{\alpha}\zeta_{\nu}\left(  y\right)
\right)  _{k}=0$.
\end{enumerate}
\end{proposition}

\begin{proof}
For part (1) suppose $\alpha_{i}\neq\beta_{i}$ for some $i$ then $\left(
x^{\alpha}g_{1}\left(  y\right)  ,x^{\beta}g_{2}\left(  y\right)  \right)
_{k}=\left(  \tau_{i}x^{\alpha}g_{1}\left(  y\right)  ,\tau_{i}x^{\beta}%
g_{2}\left(  y\right)  \right)  _{k}=\allowbreak\eta^{\beta_{i}-\alpha_{i}%
}\left(  x^{\alpha}g_{1}\left(  y\right)  ,x^{\beta}g_{2}\left(  y\right)
\right)  _{k}$ and $\eta^{\beta_{i}-\alpha_{i}}\neq1$.

For part (2), suppose that $\mu_{i}\neq\nu_{i}$ for some $i$. By Lemma \ref{lem:sadj},
we have the equality
$\left(  \mathcal{U}_{i}x^{\alpha}\zeta_{\mu}\left(  y\right)  ,x^{\alpha
}\zeta_{\nu}\left(  y\right)  \right)  _{k}=\left(  x^{\alpha}\zeta_{\mu
}\left(  y\right)  ,\mathcal{U}_{i}x^{\alpha}\zeta_{\nu}\left(  y\right)
\right)  _{k}$. This proves the claim, because the two polynomials are
eigenfunctions with different eigenvalues.
\end{proof}

For other parity types see Proposition \ref{nonstd}. To do the nontrivial
pairings we collect some facts and notation from the $S_{N}$ case. An element
$\lambda\in\mathbb{N}_{0}^{N}$ is a \textit{partition} if $\lambda_{1}%
\geq\lambda_{2}\geq\cdots\geq\lambda_{N}$. For any $\gamma\in\mathbb{N}%
_{0}^{N}$ let $\gamma^{+}$ be the unique partition such that $\gamma
^{+}=w\gamma$ for some $w\in S_{N}$. For $\gamma\in\mathbb{N}_{0}^{N}$ and
$\varepsilon=\pm$ let
\[
\mathcal{E}_{\varepsilon}\left(  \gamma\right)  =\prod\left\{  1+\frac
{\varepsilon\kappa_{0}}{\xi_{j}\left(  \gamma\right)  -\xi_{i}\left(
\gamma\right)  }:i<j\,\&\,\gamma_{i}<\gamma_{j}\right\}  .
\]
For a partition $\lambda$ and indeterminate $t$ (and implicit parameter
$\kappa_{0}$) the \textit{generalized Pochhammer symbol} is defined by
\[
\left(  t\right)  _{\lambda}=\prod_{i=1}^{N}\left(  t-\left(  i-1\right)
\kappa_{0}\right)  _{\lambda_{i}},
\]
where $\left(  a\right)  _{0}=1$ and $\left(  a\right)  _{n+1}=\left(
a\right)  _{n}\left(  a+n\right)  $; further suppose $\lambda_{M}%
>\lambda_{M+1}=0$ (or $M=N$) then the \textit{hook length product} is
\[
h\left(  \lambda,t\right)  =\prod_{i=1}^{M}\prod_{j=1}^{\lambda_{i}}%
(\lambda_{i}-j+t+\kappa_{0}\,\#\left\{  s:s>i\,\&\,j\leq\lambda_{s}\right\}
).
\]
For $0\leq l\leq m-1$ and $1\leq i\leq m-1$ define the indicator functions
$\chi_{i}\left(  l\right)  =\left\lfloor \frac{l}{i}\right\rfloor $ and extend
coordinate-wise to parity types, that is, $\chi_{i}\left(  \alpha\right)
\in\mathbb{N}_{0}^{N}$ and $\left(  \chi_{i}\left(  \alpha\right)  \right)
_{j}=\chi_{i}\left(  \alpha_{j}\right)  $ for each $j$. The rest of the
section is mostly concerned with proving the following.

\begin{theorem}
\label{bigthm}Let $\alpha$ be a standard parity type and let $\gamma
\in\mathbb{N}_{0}^{N}$, then
\begin{gather*}
\left(  x^{\alpha}\zeta_{\gamma}\left(  y\right)  ,x^{\alpha}\zeta_{\gamma
}\left(  y\right)  \right)  _{k}=m^{\left|  \alpha\right|  +m\left|
\gamma\right|  }\prod_{i=1}^{m-1}\left(  \left(  N-1\right)  \kappa_{0}%
+\frac{i}{m}+\kappa_{i}\right)  _{\left(  \gamma+\chi_{i}\left(
\alpha\right)  \right)  ^{+}}\\
\times\left(  N\kappa_{0}+1\right)  _{\gamma^{+}}\frac{h\left(  \gamma
^{+},\kappa_{0}+1\right)  }{h\left(  \gamma^{+},1\right)  }\mathcal{E}%
_{+}\left(  \gamma\right)  \mathcal{E}_{-}\left(  \gamma\right)  .
\end{gather*}
\end{theorem}

The argument has three main steps:

\begin{enumerate}
\item  find $\left(  x^{\alpha}\zeta_{\gamma}\left(  y\right)  ,x^{\alpha
}\zeta_{\gamma}\left(  y\right)  \right)  _{k}$ in terms of $\left(
\zeta_{\gamma}\left(  y\right)  ,\zeta_{\gamma}\left(  y\right)  \right)  _{k};$

\item  find $\left(  \zeta_{\gamma}\left(  y\right)  ,\zeta_{\gamma}\left(
y\right)  \right)  _{k}$ in terms of $\left(  \zeta_{\lambda}\left(  y\right)
,\zeta_{\lambda}\left(  y\right)  \right)  _{k}$, where $\lambda=\gamma^{+}$;

\item  find $\left(  \zeta_{\lambda}\left(  y\right)  ,\zeta_{\lambda}\left(
y\right)  \right)  _{k}$ in terms of $\left(  \zeta_{\lambda-\varepsilon_{M}%
}\left(  y\right)  ,\zeta_{\lambda-\varepsilon_{M}}\left(  y\right)  \right)
_{k}$ where $\lambda_{M}>0$ and $\lambda_{M+1}=0$ or $M=N. $
\end{enumerate}

The formula is valid for the trivial case $\alpha=0=\gamma$. The method of
proof is to assume the formula for lower degrees and show the above mentioned
ratios are consistent with the formula.

\begin{proposition}
Let $\alpha$ be a standard parity type and let $\gamma\in\mathbb{N}_{0}^{N}$
then
\[
T^{\alpha}x^{\alpha}\zeta_{\gamma}\left(  y\right)  =m^{\left|  \alpha\right|
}\frac{\prod_{i=1}^{m-1}\left(  \left(  N-1\right)  \kappa_{0}+\frac{i}%
{m}+\kappa_{i}\right)  _{\left(  \gamma+\chi_{i}\left(  \alpha\right)
\right)  ^{+}}}{\prod_{i=1}^{m-1}\left(  \left(  N-1\right)  \kappa_{0}%
+\frac{i}{m}+\kappa_{i}\right)  _{\gamma^{+}}}\zeta_{\gamma}\left(  y\right)
.
\]
\end{proposition}

\begin{proof}
We reduce the degree of $x^{\alpha}$ by one at each step so that each
intermediate stage is a standard parity type. Suppose that $\alpha
_{s}=l>\alpha_{s+1}$ for some $s$, or $s=N$ and $l>0$, then $\left\{  j:j\neq
s,\alpha_{j}\geq l\right\}  =\left\{  j:j<s\right\}  $. By Proposition
\ref{tixg}
\begin{align*}
T_{s}x^{\alpha}\zeta_{\gamma}\left(  y\right)   & =mx^{\alpha}x_{s}%
^{-1}\left(  \frac{l}{m}+\kappa_{l}+\mathcal{D}_{s}y_{s}-1-\kappa_{0}%
\sum_{j<s}\left(  j,s\right)  \right)  \zeta_{\gamma}\left(  y\right) \\
& =mx^{\alpha}x_{s}^{-1}\left(  \frac{l}{m}+\kappa_{l}+\xi_{s}\left(
\gamma\right)  -1-\kappa_{0}\right)  \zeta_{\gamma}\left(  y\right)  ,
\end{align*}
(because $\mathcal{D}_{s}y_{s}-\kappa_{0}\sum_{j<s}\left(  j,s\right)
=\mathcal{U}_{s}^{A}-\kappa_{0}$); recall $\xi_{s}\left(  \gamma\right)
=\kappa_{0}(N-\#\left\{  j:\gamma_{j}>\gamma_{s}\right\}  -\#\left\{
j:\gamma_{j}=\gamma_{s}\,\&\,j<s\right\}  \allowbreak+\gamma_{s}+1)$. There is
a unique $w\in S_{N}$ such that $w\gamma=\gamma^{+}$ and $i<j$ and $\gamma
_{i}=\gamma_{j}$ implies $w\left(  i\right)  <w\left(  j\right)  $. We claim
that $w\left(  \gamma+\chi_{t}\left(  \alpha\right)  \right)  =\left(
\gamma+\chi_{t}\left(  \alpha\right)  \right)  ^{+}$ for $1\leq t\leq m-1$. We
must show that $\gamma_{i}+\chi_{t}\left(  \alpha_{i}\right)  >\gamma_{j}%
+\chi_{t}\left(  \alpha_{j}\right)  $ implies $w\left(  i\right)  <w\left(
j\right)  $; the condition is equivalent to $\gamma_{i}-\gamma_{j}>\chi
_{t}\left(  \alpha_{j}\right)  -\chi_{t}\left(  \alpha_{i}\right)  $. If
$\chi_{t}\left(  \alpha_{j}\right)  \geq\chi_{t}\left(  \alpha_{i}\right)  $
then $\gamma_{i}>\gamma_{j}$ and $w\left(  i\right)  <w\left(  j\right)  $. By
definition of standard type the finite sequence $\left(  \chi_{t}\left(
\alpha_{i}\right)  \right)  _{i=1}^{N}$ is nonincreasing (a partition), thus
$\chi_{t}\left(  \alpha_{j}\right)  <\chi_{t}\left(  \alpha_{i}\right)  $
implies $\chi_{t}\left(  \alpha_{i}\right)  =1,\chi_{t}\left(  \alpha
_{j}\right)  =0$ and $i<j$. So $\gamma_{i}+1>\gamma_{j}$; if $\gamma
_{i}>\gamma_{j}$ or if $\gamma_{i}=\gamma_{j}$ and $i<j$ then $w\left(
i\right)  <w\left(  j\right)  $.

Let $\beta=\alpha-\varepsilon_{s}$ (the parity type of $T_{s}x^{\alpha}%
\zeta_{\gamma}\left(  y\right)  $), then $\chi_{t}\left(  \beta\right)
=\chi_{t}\left(  \alpha\right)  $ for $t\neq l$. By the claim the only
difference between $\left(  \gamma+\chi_{t}\left(  \alpha\right)  \right)
^{+}$ and $\left(  \gamma+\chi_{t}\left(  \beta\right)  \right)  ^{+}$ is
$\left(  \gamma+\chi_{t}\left(  \alpha\right)  \right)  _{w\left(  s\right)
}^{+}=\gamma_{s}+1=\left(  \gamma+\chi_{t}\left(  \beta\right)  \right)
_{w\left(  s\right)  }^{+}+1$ (also note that $\left(  \gamma+\chi_{t}\left(
\alpha\right)  \right)  _{w\left(  j\right)  }^{+}=\left(  \gamma+\chi
_{t}\left(  \alpha\right)  \right)  _{j}$ and similarly for $\beta$). Up to a
factor depending only on $\gamma$, and where $C_{l}$ denotes the factors
involving $\gamma+\chi_{t}\left(  \alpha\right)  $ for $t\neq l$ which do not
change from $\beta$ to $\alpha$, $T^{\beta}x^{\beta}\zeta_{\gamma}\left(
y\right)  $ equals
\begin{align*}
& C_{l}m^{\left|  \beta\right|  }\prod_{i=1}^{N}\left(  \left(  N-i\right)
\kappa_{0}+\frac{l}{m}+\kappa_{l}\right)  _{\left(  \gamma+\chi_{l}\left(
\beta\right)  \right)  _{i}^{+}}\\
& =C_{l}m^{\left|  \beta\right|  }\prod_{j=1}^{N}\left(  \left(  N-w\left(
j\right)  \right)  \kappa_{0}+\frac{l}{m}+\kappa_{l}\right)  _{\left(
\gamma+\chi_{l}\left(  \beta\right)  \right)  _{j}}.
\end{align*}
By construction $\gamma_{w\left(  i\right)  }^{+}=\gamma_{i}$ and thus
$\xi_{s}\left(  \gamma\right)  =\xi_{w\left(  s\right)  }\left(  \gamma
^{+}\right)  =\left(  N-w\left(  s\right)  +1\right)  \kappa_{0}+\gamma_{s}%
+1$. We showed that $T_{s}x^{\alpha}\zeta_{\gamma}\left(  y\right)  =m\left(
\frac{l}{m}+\kappa_{l}+\xi_{s}\left(  \gamma\right)  -1-\kappa_{0}\right)
x^{\beta}\zeta_{\gamma}\left(  y\right)  =\allowbreak m\left(  A+\gamma
_{s}\right)  x^{\beta}\zeta_{\gamma}\left(  y\right)  $, where $A=\frac{l}%
{m}+\kappa_{l}+\left(  N-w\left(  s\right)  \right)  \kappa_{0}$. This is the
desired factor since $\left(  \gamma+\chi_{l}\left(  \beta\right)  \right)
_{s}=\gamma_{s}$ and $\left(  A\right)  _{\gamma_{s}}\left(  A+\gamma
_{s}\right)  =\left(  A\right)  _{\gamma_{s}+1}$ and $\left(  \gamma+\chi
_{l}\left(  \alpha\right)  \right)  _{w\left(  s\right)  }^{+}=\left(
\gamma+\chi_{l}\left(  \alpha\right)  \right)  _{s}=\gamma_{s}+1$. This
process is used repeatedly to find $T^{\alpha}x^{\alpha}\zeta_{\gamma}\left(
y\right)  $ as a multiple of $\zeta_{\gamma}\left(  y\right)  $.
\end{proof}

We state the definition of an admissible inner product $\left\langle
\cdot,\cdot\right\rangle $ on polynomials associated with $S_{N}.$ For
polynomials $g_{1}\left(  y\right)  ,g_{2}\left(  y\right)  $ we require
that (1) $\left\langle g_{1}\left(  y\right)  ,g_{2}\left(  y\right)
\right\rangle =0$ if $g_{1},g_{2}$ are homogeneous of different degrees, (2)
for any $w\in S_{N}$,
$\left\langle wg_{1}\left(  y\right)  ,wg_{2}\left(  y\right)  \right\rangle
=\left\langle g_{1}\left(  y\right)  ,g_{2}\left(  y\right)  \right\rangle $
and (3) $\left\langle \mathcal{D}_{i}y_{i}g_{1}\left(
y\right)  ,g_{2}\left(  y\right)  \right\rangle =\left\langle g_{1}\left(
y\right)  ,\mathcal{D}_{i}y_{i}g_{2}\left(  y\right)  \right\rangle $ for
$1\leq i\leq N$. The pairing $\left\langle \cdot,\cdot\right\rangle _{T}$ is
permissible; the first two properties have been established already, and by
Proposition \ref{tixg} we have $T_{i}x_{i}g\left(  y\right)  =m\left(
\frac{1}{m}+\kappa_{1}-1+\mathcal{D}_{i}y_{i}\right)  g\left(  y\right)  $.
This proves property (3). The following is a consequence of Theorem 8.5.8 in
\cite{DX}.

\begin{proposition}
Suppose $\gamma\in\mathbb{N}_{0}^{N}$. Then we have
\[\left(  \zeta_{\gamma}\left(
y\right)  ,\zeta_{\gamma}\left(  y\right)  \right)  _{k}=\mathcal{E}%
_{+}\left(  \gamma\right)  \mathcal{E}_{-}\left(  \gamma\right)  \left(
\zeta_{\gamma^{+}}\left(  y\right)  ,\zeta_{\gamma^{+}}\left(  y\right)
\right)  _{k}.\]
\end{proposition}
Next we consider the problem of lowering the degree of $\zeta_{\lambda}$ for a
partition $\lambda$. Suppose that $\lambda_{M}>0=\lambda_{M+1}$ (where
$\lambda_{N+1}=0$; by using several results from Section 8.6 in \cite{DX} we
will compute $\dfrac{\left(  \zeta_{\lambda}\left(  y\right)  ,\zeta_{\lambda
}\left(  y\right)  \right)  _{k}}{\left(  \zeta_{\lambda-\varepsilon_{M}%
}\left(  y\right)  ,\zeta_{\lambda-\varepsilon_{M}}\left(  y\right)  \right)
_{k}}$. The starting point is the equation $\left(  T_{M}^{m}\zeta_{\lambda
}\left(  y\right)  ,T_{M}^{m}\zeta_{\lambda}\left(  y\right)  \right)
_{k}=\left(  x_{M}^{m}T_{M}^{m}\zeta_{\lambda}\left(  y\right)  ,\zeta
_{\lambda}\left(  y\right)  \right)  _{k}=\left(  y_{M}T_{M}^{m}\zeta
_{\lambda}\left(  y\right)  ,\zeta_{\lambda}\left(  y\right)  \right)  _{k}$.

We claim that $T_{M}^{m}\zeta_{\lambda}\left(  y\right)  =m^{m}K_{\lambda
}\mathcal{D}_{M}\zeta_{\lambda}\left(  y\right)  $, where
\[K_{\lambda
}=\prod_{t=1}^{m-1}\left(  \frac{t}{m}+\kappa_{t}+\left(  N-M\right)
\kappa_{0}+\lambda_{M}-1\right).\]
Indeed, $T_{M}\zeta_{\lambda}\left(
y\right)  =mx_{M}^{m-1}\mathcal{D}_{M}\zeta_{\lambda}\left(  y\right)  $. By
Lemma 8.6.3(ii) of \cite{DX} we know in addition that
$\mathcal{D}_{M}\,y_{M}\,\mathcal{D}_{M}%
\zeta_{\lambda}\left(  y\right)  =\left(  \left(  N-M\right)  \kappa
_{0}+\lambda_{M}\right)  \mathcal{D}_{M}\zeta_{\lambda}\left(  y\right)  $. By
induction suppose that
\[
T_{M}^{n}\zeta_{\lambda}\left(  y\right)  =m^{n}\prod_{t=m-n+1}^{m-1}\left(
\frac{t}{m}+\kappa_{t}+\left(  N-M\right)  \kappa_{0}+\lambda_{M}-1\right)
x_{M}^{m-n}\mathcal{D}_{M}\zeta_{\lambda}\left(  y\right).
\]
Then
\begin{align*}
T_{M}x_{M}^{m-n}&\mathcal{D}_{M}\zeta_{\lambda}\left(  y\right)
=mx_{M}^{m-n-1}\left(  \frac{m-n}{m}+\kappa_{m-n}+\mathcal{D}_{M}%
\,y_{M}-1\right)  \mathcal{D}_{M}\zeta_{\lambda}\left(  y\right) \\
& =mx_{M}^{m-n-1}\left(  \frac{m-n}{m}+\kappa_{m-n}+\left(  N-M\right)
\kappa_{0}+\lambda_{M}-1\right)  \mathcal{D}_{M}\zeta_{\lambda}\left(
y\right),
\end{align*}
where we let $n$ take the values $1,2,\ldots,m-1$ in turn.

Suppose that $\theta_{M}=\left(  1,2\right)  \left(  2,3\right)  \ldots\left(
M-1,M\right)  \in S_{N}$ is a cycle. Denote by $\widetilde{\lambda}\allowbreak =\theta
_{M}\left(  \lambda-\varepsilon_{M}\right)  =\left(  \lambda_{M}-1,\lambda
_{1},\lambda_{2},\ldots,\lambda_{M-1},0,\ldots\right)  $. By Theorem 8.6.5 of
\cite{DX} $\mathcal{D}_{M}\zeta_{\lambda}=b_{\lambda}\theta_{M}^{-1}%
\zeta_{\widetilde{\lambda}}$, where $b_{\lambda}=\left(  N-M+1\right)
\kappa_{0}+\lambda_{M}$, and by Lemma 8.6.2(iii) $\left(  y_{M}\mathcal{D}%
_{M}\zeta_{\lambda},\zeta_{\lambda}\right)  _{k}=\dfrac{\lambda_{M}}%
{\kappa_{0}+\lambda_{M}}b_{\lambda}\left(  \zeta_{\lambda},\zeta_{\lambda
}\right)  _{k}$. Combining these equations we obtain
\begin{align*}
\left(  T_{M}^{m}\zeta_{\lambda}\left(  y\right)  ,T_{M}^{m}\zeta_{\lambda
}\left(  y\right)  \right)  _{k}  & =m^{2m}K_{\lambda}^{2}\left(
\mathcal{D}_{M}\zeta_{\lambda},\mathcal{D}_{M}\zeta_{\lambda}\right)
_{k}=\left(  y_{M}T_{M}^{m}\zeta_{\lambda}\left(  y\right)  ,\zeta_{\lambda
}\left(  y\right)  \right)  _{k}\\
& =m^{m}K_{\lambda}\left(  y_{M}\,\mathcal{D}_{M}\zeta_{\lambda}%
,\zeta_{\lambda}\right)  _{k}=m^{m}K_{\lambda}b_{\lambda}\dfrac{\lambda_{M}%
}{\kappa_{0}+\lambda_{M}}\left(  \zeta_{\lambda},\zeta_{\lambda}\right)  _{k},
\end{align*}
and thus
\begin{align*}
\left(  \zeta_{\lambda},\zeta_{\lambda}\right)  _{k}  & =m^{m}K_{\lambda
}b_{\lambda}\dfrac{\kappa_{0}+\lambda_{M}}{\lambda_{M}}\left(  \theta_{M}%
^{-1}\zeta_{\widetilde{\lambda}},\theta_{M}^{-1}\zeta_{\widetilde{\lambda}%
}\right)  _{k}\\
& =m^{m}K_{\lambda}b_{\lambda}\dfrac{\kappa_{0}+\lambda_{M}}{\lambda_{M}%
}\left(  \zeta_{\lambda-\varepsilon_{M}},\zeta_{\lambda-\varepsilon_{M}%
}\right)  _{k}\,\mathcal{E}_{+}\left(  \widetilde{\lambda}\right)
\mathcal{E}_{-}\left(  \widetilde{\lambda}\right)  ,
\end{align*}
because $\left(  \theta_{M}^{-1}\zeta_{\widetilde{\lambda}},\theta_{M}%
^{-1}\zeta_{\widetilde{\lambda}}\right)  _{k}=\left(  \zeta_{\widetilde
{\lambda}},\zeta_{\widetilde{\lambda}}\right)  _{k}$ and $\left(
\widetilde{\lambda}\right)  ^{+}=\lambda-\varepsilon_{M}$. Now notice that $\mathcal{E}%
_{\varepsilon}\left(  \widetilde{\lambda}\right)  =\prod_{i=2}^{M}\left(
1+\dfrac{\varepsilon\kappa_{0}}{\left(  M-i\right)  \kappa_{0}+\lambda
_{i}+1-\lambda_{M}}\right)  $ for $\varepsilon=\pm$, and it can be shown (see
Section 8.7 of \cite{DX}) that
\[
\dfrac{\kappa_{0}+\lambda_{M}}{\lambda_{M}}\mathcal{E}_{+}\left(
\widetilde{\lambda}\right)  \mathcal{E}_{-}\left(  \widetilde{\lambda}\right)
=\frac{h\left(  \lambda,\kappa_{0}+1\right)  h\left(  \lambda-\varepsilon
_{M},1\right)  }{h\left(  \lambda-\varepsilon_{M},\kappa_{0}+1\right)
h\left(  \lambda,1\right)  }.
\]
Further $K_{\lambda}=\prod_{t=1}^{m-1}\dfrac{\left(  \left(  N-1\right)
\kappa_{0}+\frac{t}{m}+\kappa_{t}\right)  _{\lambda}}{\left(  \left(
N-1\right)  \kappa_{0}+\frac{t}{m}+\kappa_{t}\right)  _{\lambda-\varepsilon
_{M}}}$ and $b_{\lambda}=\dfrac{\left(  N\kappa_{0}+1\right)  _{\lambda}%
}{\left(  N\kappa_{0}+1\right)  _{\lambda-\varepsilon_{M}}}$ . This completes
the proof of Theorem \ref{bigthm}.

We turn to skew-symmetric polynomials in $y$. Suppose $\lambda$ is a partition
with all distinct parts, that is, $\lambda_{1}>\lambda_{2}>\ldots>\lambda
_{N}\geq0$ then the polynomial $a_{\lambda}=\sum_{w\in S_{N}}\mathrm{sign}%
\left(  w\right)  w\zeta_{\lambda}=\sum_{w\in S_{N}}\mathrm{sign}\left(
w\right)  \dfrac{\mathcal{E}_{-}\left(  \lambda^{R}\right)  }{\mathcal{E}%
_{-}\left(  w\lambda\right)  }\zeta_{w\lambda}$ is skew-symmetric (that is,
$\left(  i,j\right)  a_{\lambda}=-a_{\lambda}$ for any transposition), where
$\lambda^{R}=\left(  \lambda_{N},\lambda_{N-1},\ldots,\lambda_{1}\right)  $;
further $\left(  a_{\lambda},a_{\lambda}\right)  _{k}=N!\mathcal{E}_{-}\left(
\lambda^{R}\right)  \left(  \zeta_{\lambda},\zeta_{\lambda}\right)  _{k}$ (see
Theorem 8.5.11 in \cite{DX}). For the minimal example $\lambda=\delta=\left(
N-1,N-2,\ldots,2,1,0\right)  $ the polynomial $a_{\delta}$ has the simple form
$\prod_{1\leq i<j\leq N}\left(  y_{i}-y_{j}\right)  $ and $\mathcal{E}%
_{-}\left(  \delta^{R}\right)  =\dfrac{h\left(  \delta,1\right)  }{h\left(
\delta,\kappa_{0}+1\right)  }$. Let $\upsilon=\left(  1,1,\ldots,1\right)
\in\mathbb{N}_{0}^{N}$ and for any $0\leq t\leq m-1$ let $f\left(  x\right)
=x^{t\upsilon}\prod_{1\leq i<j\leq N}\left(  x_{i}^{m}-x_{j}^{m}\right)  $ (so
that the parity type is $\left(  t,t,\ldots,t\right)  .$ Specializing Theorem
\ref{bigthm} to $\gamma=\delta$ and $\alpha=t\upsilon$, note that $a_{\delta
}\left(  y\right)  =\sum_{w\in S_{N}}\mathrm{sign}\left(  w\right)
w\zeta_{\delta}\left(  y\right)  $ and thus
\begin{align*}
T^{t\upsilon}x^{t\upsilon}a_{\delta}\left(  y\right)   & =\sum_{w\in S_{N}%
}\mathrm{sign}\left(  w\right)  w\left(  T^{t\upsilon}x^{t\upsilon}%
\zeta_{\delta}\left(  y\right)  \right) \\
& =m^{mt}\prod_{i=1}^{t}\frac{\left(  \left(  N-1\right)  \kappa_{0}+\frac
{i}{m}+\kappa_{i}\right)  _{\delta+\upsilon}}{\left(  \left(  N-1\right)
\kappa_{0}+\frac{i}{m}+\kappa_{i}\right)  _{\delta}}a_{\delta}\left(
y\right)
\end{align*}
which implies
\begin{align*}
\left(  f,f\right)  _{k}  & =N!\,m^{m\left(  t+N\left(  N-1\right)  /2\right)
}\prod_{i=1}^{t}\left(  \left(  N-1\right)  \kappa_{0}+\frac{i}{m}+\kappa
_{i}\right)  _{\delta+\upsilon}\\
& \times\prod_{i=t+1}^{m-1}\left(  \left(  N-1\right)  \kappa_{0}+\frac{i}%
{m}+\kappa_{i}\right)  _{\delta}\left(  N\kappa_{0}+1\right)  _{\delta}.
\end{align*}
We used the fact that $\left\langle f,g\right\rangle /\left\langle
\zeta_{\lambda},\zeta_{\lambda}\right\rangle $ is the same in any permissible
inner product for each $f,g\in\mathrm{span}\left\{  \zeta_{w\lambda}:w\in
S_{N}\right\}  $ and the value of $\left\langle a_{\delta},a_{\delta
}\right\rangle /\left\langle \zeta_{\delta},\zeta_{\delta}\right\rangle $ from
Theorem 8.7.15 in \cite{DX}. The formula was conjectured by P. Hanlon in 1995.

\subsection{The Radical}

The radical for the hermitian form $(\cdot,\cdot)_k$
is the linear space $\mathrm{Rad}\left(
\kappa\right)  $ of polynomials $p$ such that $\left(  q,p\right)  _{k}=0$ for
any polynomial $q$; the space depends on the parameter values $\kappa
_{0},\kappa_{1},\ldots,\kappa_{m-1}$. For the group $G\left(  m,p,N\right)  $
the set of singular values (when $\mathrm{Rad}\left(  \kappa\right)
\neq\left\{  0\right\}  $) can be explicitly stated with reference to the
discriminant. Briefly, the radical is nontrivial exactly when any of the
linear functions occurring in $\left(  f,f\right)  _{k}$ take on values in
$0,-1,-2,-3\ldots$ where $f=\left(  x_{1}x_{2}\ldots x_{N}\right)  ^{m-1}%
\prod_{1\leq i<j\leq N}\left(  x_{i}^{m}-x_{j}^{m}\right)  $.
The proof of
this requires extra care when $\kappa_{0}\in-\mathbb{N}$ because some of the
nonsymmetric Jack polynomials fail to exist for such values. By careful
analysis of their construction one sees that the poles can occur only at
values of the form $n\kappa_{0}+l=0$ for $1\leq n\leq N$ and $l\in\mathbb{N}$.
In the following $\kappa$ refers to the parameters $\left(  \kappa_{i}\right)
_{i=0}^{m-1}$

\begin{definition}
Let $\mathcal{K}_{0}=\left\{  \kappa:\kappa_{0}=-\frac{j}{n}-l:l\in
\mathbb{N}_{0},2\leq j-1\leq n\leq N\right\}  $ and let $\mathcal{K}_{1}%
=\cup_{i=1}^{m-1}\cup_{n=0}^{N-1}\left\{  \kappa:n\kappa_{0}+\frac{i}%
{m}+\kappa_{i}\in-\mathbb{N}_{0}\right\}  $.
\end{definition}

We claim the set of singular values is $\mathcal{K}_{0}\cup\mathcal{K}_{1}$.

Firstly, suppose $\kappa\notin\mathcal{K}_{0}$ and $\kappa_{0}\notin
-\mathbb{N}$; in this case all the functions $\zeta_{\alpha}$ exist, the
simultaneous eigenfunctions of $\left\{  \mathcal{U}_{i}\right\}  _{i=1}^{N}$
span the polynomials and the pairing formulae are valid. If $\kappa
\in\mathcal{K}_{1}$ then there is a (nonzero) polynomial $x^{\alpha}%
\zeta_{\lambda}\left(  y\right)  $ (some standard parity type $\alpha$, a
partition $\lambda$) which is orthogonal to each polynomial, thus in
$\mathrm{Rad}\left(  \kappa\right)  $. On the other hand, suppose there is a
nonzero $p\in\mathrm{Rad}\left(  \kappa\right)  $ then expand $p$ in terms of
the eigenvector basis and let $wx^{\alpha}\zeta_{\gamma}\left(  y\right)  $
appear in the expansion with a nonzero coefficient (see Proposition
\ref{nonstd}); by the orthogonality relations this implies $\left(  x^{\alpha
}\zeta_{\gamma}\left(  y\right)  ,x^{\alpha}\zeta_{\gamma}\left(  y\right)
\right)  _{k}=0$ and thus $\kappa\in\mathcal{K}_{1}$.

Secondly, suppose $\kappa\in\mathcal{K}_{1}$ and $\kappa_{0}=l\in-\mathbb{N}$;
but $\mathcal{K}_{1}$ is closed so take a sufficiently close $\kappa^{\prime
}\in\mathcal{K}_{1}$ with $\kappa_{0}^{\prime}\notin-\mathbb{N}$ (in the same
component as $\kappa$), for this value there is a nonzero polynomial $p$ (a
simultaneous eigenfunction) in $\mathrm{Rad}\left(  \kappa\right)  $; for some
$n=1,2,3,\ldots$ $\left(  \kappa_{0}-l\right)  ^{n}p$ has no poles and a
nonzero limit at $\kappa_{0}=l$. For any polynomial $q$ we have $\left(
q,\left(  \kappa_{0}-l\right)  ^{n}p\right)  _{k}=0$ for all $\kappa^{\prime}$
close to $\kappa$, thus $\mathrm{Rad}\left(  \kappa\right)  $ is nontrivial.

Thirdly, let $\kappa\in\mathcal{K}_{0}$; by the results in \cite{DjO} there is
a nonzero polynomial $g\left(  y\right)  $ of least degree in the type-$A$
radical, which implies $\mathcal{D}_{i}g\left(  y\right)  =0$ for $1\leq i\leq
N$. By Proposition \ref{tixg} $T_{i}g\left(  y\right)  =0$ for all $i$ and so
$g\left(  y\right)  \in\mathrm{Rad}\left(  \kappa\right)  .$

It remains to show that if $\kappa\notin\mathcal{K}_{0}\cup\mathcal{K}_{1}$
and $\kappa_{0}\in-\mathbb{N}$ then $\mathrm{Rad}\left(  \kappa\right)
=\left\{  0\right\}  $. The problem is that the simultaneous eigenfunctions of
$\left\{  \mathcal{U}_{i}\right\}  _{i=1}^{N}$ no longer span the polynomials;
nevertheless it is still possible to give an argument based on triangularity
properties. We will show that $\kappa\notin\mathcal{K}_{1}$ and $\mathrm{Rad}%
\left(  \kappa\right)  \neq\left\{  0\right\}  $ implies $\kappa\in
\mathcal{K}_{0}$. Suppose that $p\in\mathrm{Rad}\left(  \kappa\right)  $ and
$p\neq0$; because $\mathrm{Rad}\left(  \kappa\right)  $ is an ideal the
polynomial $g\left(  y\right)  =x^{\alpha}p\left(  x\right)  \in
\mathrm{Rad}\left(  \kappa\right)  $ (if the parity type of $p$ is $\beta$ let
$\alpha_{i}=m-\beta_{i}$ for each $i$). Suppose that $g_{0}\left(  y\right)
\in\mathrm{Rad}\left(  \kappa\right)  $ and has minimal degree (among
polynomials in $y$). Consider the polynomials $T_{i}g_{0}\left(  y\right)
=mx_{i}^{m-1}\mathcal{D}_{i}g_{0}\left(  y\right)  $; if each $\mathcal{D}%
_{i}g_{0}=0$ then $g_{0}$ is in the type-$A$ radical and $\kappa\in
\mathcal{K}_{0}$. Suppose not, we may assume $\mathcal{D}_{1}g_{0}\left(
y\right)  \neq0$. By Proposition \ref{tixg} $T_{1}^{m}g_{0}\left(  y\right)
=m^{m}\prod_{i=1}^{m-1}\left(  \frac{i}{m}+\kappa_{i}+\mathcal{D}_{1}%
y_{1}-1\right)  \mathcal{D}_{1}g_{0}\left(  y\right)  $. We now use the fact
that $\mathcal{D}_{1}y_{1}$ is triangular (see p.454 in \cite{D5}) with
respect to the partial order $\vartriangleright$ (see Definition
\ref{defdom}). Let $cy^{\gamma}$ be a nonzero term in $\mathcal{D}_{1}g_{0}$
with maximal (for $\vartriangleright$) $\gamma;$this implies that the
coefficient of $y^{\gamma}$ in $T_{1}^{m}g_{0}\left(  y\right)  $ is
$m^{m}\prod_{i=1}^{m-1}\left(  \frac{i}{m}+\kappa_{i}+v\right)  c$ where
$\nu=\kappa_{0}\left(  N-1-\#\left\{  j:\gamma_{j}>\gamma_{1}\right\}
\right)  +\gamma_{1}$; the hypothesis $\kappa\notin\mathcal{K}_{1}$ implies
$\frac{i}{m}+\kappa_{i}+v\neq0$ for each $i$ but then $T_{1}^{m}g_{0}\left(
y\right)  $ is a nonzero polynomial in $y$ in the radical and of lower degree
than $g_{0}$, a contradiction.

Thus the detailed knowledge of type-$A$ polynomials makes it possible to
describe the set of singular values for $G\left(  m,1,N\right)  $ and indeed
for any $G\left(  m,p,N\right)  $. For the latter impose the periodicity
conditions in Remark \ref{GmpN} on $\kappa$.

\subsection{Shift Operators}

Suppose $G$ is a complex reflection group and recall the definitions of
Section 2. Let $p$ be a $G$-invariant
polynomial, and $C$ a $G$-orbit of reflection hyperplanes. Given a rational
function $f$ on $V$, we denote by $m(f)$ the operator ``multiplication by $f$''.
Observe that, for
$s\in\mathbb{N}$ with $s<e_{{C}}$,  the operator
\[m\left(\prod_{H\in
{C}}\alpha_{H} ^{-s}\right) \circ p^{\ast}\left(  T\left(  k\right)
\right) \circ m\left(  \prod_{H\in{C}}\alpha_{H}^{s}\right)  \]
has the property that it maps $P^G$ to $P^G$.
\begin{qu} Is the above operator on $P^G$ equal to the restriction to $P^G$
of an operator of the form $p^{\ast}\left(  T\left(  k^{\prime}\right)  \right)$,
where $k^{\prime}$ is obtained by incrementing some
of the values of $k=\left(  k_{{C},i}\right)  $?
\end{qu}
This is well known in the case of real reflection groups \cite{H}, and in that
case the relation plays an important role in many applications \cite{O1}.
The proof of this ``shift relation'' is based in this case on the presence
of the invariant $p=\sum x_i^2$ of order two. For
this invariant the relation can be checked by simple direct computation.
Then one remarks that this forces the relation also to be true for the
higher order invariants, using $\mathfrak{sl}_2$ representation theory
(see \cite{H}).

We do not know of any general argument that works in the present case of
complex reflection groups.
Nevertheless, in this section we shall show that the answer to this question is
affirmative for the groups $G\left(  m,p,N\right)$. The argument is again based
on a reduction to the case of $S_N$.

First we deal with shifting the parameters $\kappa_{i},1\leq i\leq m-1,$ for
the group $G\left(  m,1,N\right)  $. We recall the type-$A$ commutation
$\mathcal{D}_{i}y_{i}-y_{i}\mathcal{D}_{i}=1+\kappa_{0}\sum_{j\neq i}\left(
i,j\right)  $, where $y=\left(  x_{1}^{m},\ldots,x_{N}^{m}\right)  $; thus
$y_{i}\mathcal{D}_{i}=Y_{i}$ with $Y_{i}=\mathcal{D}_{i}y_{i}-1-\kappa_{0}%
\sum_{j\neq i}\left(  i,j\right)  $ and $\left(  \mathcal{D}_{i}y_{i}\right)
\mathcal{D}_{i}=\mathcal{D}_{i}Y_{i}$, for each $i$. Let $v=\left(
1,1,\ldots,1\right)  \in\mathbb{N}^{N}$.

\begin{proposition}
Let $1\leq t\leq m-1$ and let $g$ be any polynomial in $y$ then
\[
T_{i}\left(  \kappa\right)  ^{m}x^{tv}g\left(  y\right)  =x^{tv}T_{i}\left(
\kappa^{\prime}\right)  ^{m}g\left(  y\right)  ,
\]
for $1\leq i\leq N$, where $\kappa^{\prime}=\left(  \kappa_{0},\kappa
_{1}+1,\ldots,\kappa_{t}+1,\kappa_{t+1},\ldots,\kappa_{m-1}\right)  .$
\end{proposition}

\begin{proof}
By Proposition \ref{tixg}
\begin{align*}
&  T_{i}\left(  \kappa\right)  ^{m}x^{tv}g\left(  y\right) \\
&  =m^{m}x^{tv}\prod_{s=m-t}^{m-1}\left(  \frac{s}{m}+\kappa_{s}%
+\mathcal{D}_{i}y_{i}-1\right)  \mathcal{D}_{i}\prod_{s=1}^{t}\left(  \frac
{s}{m}+\kappa_{s}+Y_{i}\right) \\
&  =m^{m}x^{tv}\mathcal{D}_{i}\prod_{s=m-t}^{m-1}\left(  \frac{s}{m}%
+\kappa_{s}+Y_{i}-1\right)  \prod_{s=1}^{t}\left(  \frac{s}{m}+\kappa
_{s}+Y_{i}\right) \\
&  =m^{m}x^{tv}\mathcal{D}_{i}\prod_{s=1}^{m-1}\left(  \frac{s}{m}+\kappa
_{s}^{\prime}+Y_{i}-1\right) \\
&  =x^{tv}T_{i}\left(  \kappa^{\prime}\right)  ^{m}g\left(  y\right)  .
\end{align*}
In each of the products the terms commute pairwise so the order does not
matter. The relation $\left(  \mathcal{D}_{i}y_{i}\right)  \mathcal{D}%
_{i}=\mathcal{D}_{i}Y_{i}$ is used to move $\mathcal{D}_{i}$ to the front
(last in the order of operation). The formula for $T_{i}\left(  \kappa
^{\prime}\right)  ^{m}g\left(  y\right)  $ is obtained by setting $t=0$ in the
starting calculation.
\end{proof}

The argument must be modified for $G\left(  m,p,N\right)  $ for $1<p<m$.
Recall from Remark \ref{GmpN} that $\kappa_{sm/p}=0$ and $\kappa
_{t+sm/p}=\kappa_{t}$ for $1\leq s\leq p-1$ and $1\leq t\leq\frac{m}{p}-1$. We
will use $\kappa$ to denote $\left(  \kappa_{0},\kappa_{1},\ldots
,\kappa_{m/p-1}\right)  $ in this discussion. The invariants for $G\left(
m,p,N\right)  $ are generated by the elementary symmetric functions of degrees
$1,2,\ldots,N-1$ in $y$ and $x^{\left(  m/p\right)  v}$. Here is the
modification of the previous proposition.

\begin{proposition}\label{prop:cycl}
For the group $G\left(  m,p,N\right)  $ let $1\leq t\leq\frac{m}{p}-1,\,0\leq
s\leq p-1$ and let $g$ be any polynomial in $y$ then
\[
T_{i}\left(  \kappa\right)  ^{m/p}x^{\left(  t+sm/p\right)  v}g\left(
y\right)  =x^{tv}T_{i}\left(  \kappa^{\prime}\right)  ^{m/p}x^{\left(
sm/p\right)  v}g\left(  y\right)  ,
\]
where $\kappa^{\prime}=\left(  \kappa_{0},\kappa_{1}+\frac{1}{p},\ldots
,\kappa_{t}+\frac{1}{p},\kappa_{t+1},\ldots,\kappa_{m/p-1}\right)  $.
\end{proposition}

\begin{proof}
Suppose first that $1\leq s\leq p-1$ then
\begin{align*}
&  T_{i}\left(  \kappa\right)  ^{m/p}x^{\left(  t+sm/p\right)  v}g\left(
y\right) \\
&  =m^{m/p}x^{\left(  t+\left(  s-1\right)  m/p\right)  v}\prod_{j=t+1+\left(
s-1\right)  m/p}^{t+sm/p}\left(  \frac{j}{m}+\kappa_{j}+Y_{i}\right)  g\left(
y\right) \\
&  =m^{m/p}x^{\left(  t+\left(  s-1\right)  m/p\right)  v}\left(  \frac{s}%
{p}+Y_{i}\right)  \prod_{j=t+1}^{m/p-1}\left(  \frac{s-1}{p}+\frac{j}%
{m}+\kappa_{j}+Y_{i}\right) \\
&  \times\prod_{j=1}^{t}\left(  \frac{s}{p}+\frac{j}{m}+\kappa_{j}%
+Y_{i}\right)  g\left(  y\right) \\
&  =m^{m/p}x^{\left(  t+\left(  s-1\right)  m/p\right)  v}\left(  \frac{s}%
{p}+Y_{i}\right)  \prod_{j=1}^{m/p-1}\left(  \frac{s-1}{p}+\frac{j}{m}%
+\kappa_{j}^{\prime}+Y_{i}\right)  g\left(  y\right) \\
&  =x^{tv}T_{i}\left(  \kappa^{\prime}\right)  ^{m/p}x^{\left(  sm/p\right)
v}g\left(  y\right)  .
\end{align*}
Again the order of the product does not matter because each term is a linear
function of the operator $Y_{i}$. The periodicity conditions were applied to
yield the third line.

Further (the case $s=0$)
\begin{align*}
&  T_{i}\left(  \kappa\right)  ^{m/p}x^{tv}g\left(  y\right) \\
&  =m^{m/p}x^{tv}\prod_{j=m-m/p+t+1}^{m-1}\left(  \frac{j}{m}+\kappa
_{j}+\mathcal{D}_{i}y_{i}-1\right)  \mathcal{D}_{i}\prod_{j=1}^{t}\left(
\frac{j}{m}+\kappa_{j}+Y_{i}\right)  g\left(  y\right) \\
&  =m^{m/p}x^{tv}\mathcal{D}_{i}\prod_{j=t+1}^{m/p-1}\left(  \frac{p-1}%
{p}+\frac{j}{m}+\kappa_{j}+Y_{i}-1\right) \\
&  \times\prod_{j=1}^{t}\left(  \frac{j}{m}+\kappa_{j}+Y_{i}\right)  g\left(
y\right) \\
&  =m^{m/p}x^{tv}\mathcal{D}_{i}\prod_{j=1}^{m/p-1}\left(  \frac{j}{m}%
-\frac{1}{p}+\kappa_{j}^{\prime}+Y_{i}\right)  g\left(  y\right) \\
&  =x^{tv}T_{i}\left(  \kappa^{\prime}\right)  ^{m/p}g\left(  y\right)  .
\end{align*}
In the second line the product over $m-\frac{m}{p}+t+1\leq j\leq m-1$ is
changed to $t+1\leq j\leq\frac{m}{p}-1$ by replacing $j$ by $j+\frac{m\left(
p-1\right)  }{p}$ and the periodicity of $\kappa_{j}$ is also used. This
completes the proof.
\end{proof}
\begin{rem}\label{rem:conv}
It may appear that the shift of $\frac{1}{p}$ is significantly different from
the shift of $1$ for the group $G\left(  m,1,N\right)  $, but in fact, the
parameters $\kappa_{i}$ should be replaced by $p\kappa_{i}$ (for $1\leq
i\leq\frac{m}{p}-1$) to conform to the general setup of Section 2 and so the
shift is effectively $1$.
\end{rem}
We turn to the parameter $\kappa_{0}$ associated with the action of the
symmetric group. Using the known results for $S_{N}$,  we are going to prove that $f\left(
T\left(  \kappa\right)  \right)  a_{\delta}\left(  y\right)  g\left(
y\right)\allowbreak  =a_{\delta}\left(  y\right)  f\left(  T\left(  \kappa^{\prime
}\right)  \right)  g\left(  y\right)  $, where $\delta=\left(  N-1,N-2,\ldots
,1,0\right)  \in\mathbb{N}_{0}^{N}$, and $a_{\delta}\left(  y\right)
=\allowbreak\prod_{1\leq i<j\leq N}\left(  y_{i}-y_{j}\right)  $,
$\kappa^{\prime}=\left(  \kappa_{0}+1,\kappa_{1},\ldots,\kappa_{m-1}\right)  $
and $f,g$ are (real) symmetric polynomials in $y$.

Recall some facts from Section 3.3: for any partition $\lambda\in
\mathbb{N}_{0}^{N}$ the space $X_{\lambda}\left(  \kappa_{0}\right)
=\mathrm{span}\left\{  w\zeta_{\lambda}:w\in S_{N}\right\}  =\mathrm{span}%
\left\{  \zeta_{w\lambda}:w\in S_{N}\right\}  $ (the $\mathbb{R}$-span) is
equipped with two permissible inner products $\left(  \cdot,\cdot\right)
_{k}$ and $\left\langle \cdot,\cdot\right\rangle _{\kappa_{0}}$, where the
inner product
$\left\langle f,g\right\rangle _{\kappa_{0}}=f\left(  \mathcal{D}\left(
\kappa_{0}\right)  \right)  g\left(  y\right)  $ for $f,g\in X_{\lambda
}\left(  \kappa_{0}\right)  $. We have shown that%
\[
\left(  f,g\right)  _{k}=m^{m|\lambda|}\prod_{i=1}^{m-1}\left(  \left(
N-1\right)  \kappa_{0}+\frac{i}{m}+\kappa_{i}\right)  _{\lambda}\left\langle
f,g\right\rangle _{\kappa_{0}}.
\]
First we establish the formula $\left(  a_{\delta}f,a_{\delta}g\right)
_{k}=\left(  a_{\delta},a_{\delta}\right)  _{k}\left(  f,g\right)
_{k^{\prime}}$ where $\left(  f,g\right)  _{k^{\prime}}$ denotes the pairing
for $\kappa^{\prime}$ and $f,g$ are symmetric polynomials in $y$. The result
of Heckman \cite{H} shows that $\left\langle a_{\delta}f,a_{\delta
}g\right\rangle _{\kappa_{0}}=\left\langle a_{\delta},a_{\delta}\right\rangle
_{\kappa_{0}}\left\langle f,g\right\rangle _{\kappa_{0}+1}$. Further, there is
a unique $S_{N}$-invariant $j_{\lambda}$ (up to scalar multiple) in
$X_{\lambda}$, and if $\lambda_{1}>\lambda_{2}>\ldots>\lambda_{N}$ then there
is a unique skew-invariant $a_{\lambda}\in X_{\lambda}$ (see Section 3.3).
Further the invariant polynomial $a_{\lambda+\delta}\left(  \kappa_{0}\right)
/a_{\delta}=$ $j_{\lambda}\left(  \kappa_{0}+1\right)  $ (note the dependence
on the parameter), see Opdam \cite{O4}; we take the constant as $1$ for
convenience. If $\lambda,\mu$ are partitions in $\mathbb{N}_{0}^{N}$ then
$\lambda\neq\mu$ implies $X_{\lambda}\left(  \kappa_{0}\right)  \bot X_{\mu
}\left(  \kappa_{0}\right)  $ in any permissible inner product, hence $\left(
a_{\lambda+\delta}\left(  \kappa_{0}\right)  ,a_{\mu+\delta}\left(  \kappa
_{0}\right)  \right)  _{k}=0$. The polynomials $a_{\lambda+\delta}\left(
\kappa_{0}\right)  /a_{\delta}$ form a basis for the symmetric polynomials so
it suffices to prove the formula for the cases $f=g=a_{\lambda+\delta}\left(
\kappa_{0}\right)  /a_{\delta}$. Thus%
\begin{align*}
&  \left(  a_{\lambda+\delta}\left(  \kappa_{0}\right)  ,a_{\lambda+\delta
}\left(  \kappa_{0}\right)  \right)  _{k}\\
&  =m^{m|\lambda+\delta|}\prod_{i=1}^{m-1}\left(  \left(  N-1\right)
\kappa_{0}+\frac{i}{m}+\kappa_{i}\right)  _{\lambda+\delta}\left\langle
a_{\lambda+\delta}\left(  \kappa_{0}\right)  ,a_{\lambda+\delta}\left(
\kappa_{0}\right)  \right\rangle _{\kappa_{0}}\\
&  =m^{m|\lambda+\delta|}\prod_{i=1}^{m-1}\prod_{j=1}^{N}\left(  \left(
N-j\right)  \kappa_{0}+\frac{i}{m}+\kappa_{i}\right)  _{\lambda_{j}+N-j}\\
&  \times\left\langle a_{\delta},a_{\delta}\right\rangle _{\kappa_{0}%
}\left\langle j_{\lambda}\left(  \kappa_{0}+1\right)  ,j_{\lambda}\left(
\kappa_{0}+1\right)  \right\rangle _{\kappa_{0}+1}\\
&  =m^{m|\lambda+\delta|}\prod_{i=1}^{m-1}\prod_{j=1}^{N}\left(  \left(
N-j\right)  \kappa_{0}+\frac{i}{m}+\kappa_{i}\right)  _{N-j}\left(  \left(
N-j\right)  \left(  \kappa_{0}+1\right)  +\frac{i}{m}+\kappa_{i}\right)
_{\lambda_{j}}\\
&  \times\left\langle a_{\delta},a_{\delta}\right\rangle _{\kappa_{0}%
}\left\langle j_{\lambda}\left(  \kappa_{0}+1\right)  ,j_{\lambda}\left(
\kappa_{0}+1\right)  \right\rangle _{\kappa_{0}+1}\\
&  =m^{m|\delta|}\prod_{i=1}^{m-1}\left(  \left(  N-1\right)  \kappa_{0}%
+\frac{i}{m}+\kappa_{i}\right)  _{\delta}\left\langle a_{\delta},a_{\delta
}\right\rangle _{\kappa_{0}}\left(  j_{\lambda}\left(  \kappa_{0}+1\right)
,j_{\lambda}\left(  \kappa_{0}+1\right)  \right)  _{k^{\prime}}%
\end{align*}
which used the relation between $\kappa_{0}+1$ and $k^{\prime}$ inner products
on $X_{\lambda}\left(  \kappa_{0}+1\right)  $, and the known $S_{N}$ results.
Setting $\lambda=0$ shows again that \newline $\left(  a_{\delta},a_{\delta
}\right)  _{k}=m^{m|\delta|}\prod_{i=1}^{m-1}\left(  \left(  N-1\right)
\kappa_{0}+\frac{i}{m}+\kappa_{i}\right)  _{\delta}\left\langle a_{\delta
},a_{\delta}\right\rangle _{\kappa_{0}}$.

\begin{proposition}\label{prop:sn}
Suppose $f,g$ are symmetric polynomials in $y$ then%
\[
f\left(  T\left(  \kappa\right)  \right)  a_{\delta}\left(  y\right)  g\left(
y\right)  =a_{\delta}\left(  y\right)  f\left(  T\left(  \kappa^{\prime
}\right)  \right)  g\left(  y\right)  ,
\]
where $\kappa^{\prime}=\left(  \kappa_{0}+1,\kappa_{1},\ldots,\kappa
_{m-1}\right)  $.
\end{proposition}

\begin{proof}
Without loss of generality assume that $f,g$ are homogeneous. By the
group-invariance properties of the pairing $f\left(  T\left(  \kappa\right)
\right)  a_{\delta}\left(  y\right)  g\left(  y\right)  $ is a skew-~symmetric
polynomial and hence divisible by $a_{\delta}\left(  y\right)  $ with a
symmetric quotient. If $\deg f>\deg g$ the quotient is zero. If
$\deg f=\deg g$ then $f\left(  T\left(  \kappa\right)  \right)  a_{\delta
}\left(  y\right)  g\left(  y\right)  =ca_{\delta}\left(  y\right)  $ for some
constant $c$; thus by the preceding formula we obtain
$a_{\delta}\left(  T\left(
\kappa\right)  \right)  f\left(  T\left(  \kappa\right)  \right)  a_{\delta
}\left(  y\right)  g\left(  y\right)  =\left(  a_{\delta}\left(  T\left(
\kappa\right)  \right)  a_{\delta}\left(  y\right)  \right)  \left(  f\left(
T\left(  \kappa^{\prime}\right)  \right)  g\left(  y\right)  \right)  $
(recall $\left(  a_{\delta},a_{\delta}\right)  _{k}\neq0$ for generic $\kappa
$, which suffices to prove this $\mathbb{Q}\left[  \kappa\right]  $ polynomial identity).

Suppose $\deg f=l<\deg g=n$ and let $f_{1}$ be an arbitrary homogeneous
symmetric polynomial of degree $n-l$.  Then we have $f_{1}\left(  T\left(
\kappa\right)  \right)  f\left(  T\left(  \kappa\right)  \right)  a_{\delta
}\left(  y\right)  g\left(  y\right)\allowbreak  =a_{\delta}\left(  y\right)
f_{1}\left(  T\left(  \kappa^{\prime}\right)  \right)  f\left(  T\left(
\kappa^{\prime}\right)  \right)  g\left(  y\right)  $, but $f\left(  T\left(
\kappa^{\prime}\right)  \right)  g\left(  y\right)  $ is symmetric of degree
$n-l$. Hence $a_{\delta}\left(  y\right)  f_{1}\left(  T\left(  \kappa^{\prime
}\right)  \right)  f\left(  T\left(  \kappa^{\prime}\right)  \right)  g\left(
y\right)  =f_{1}\left(  T\left(  \kappa\right)  \right)  a_{\delta}\left(
y\right)  f\left(  T\left(  \kappa^{\prime}\right)  \right)  g\left(
y\right)  $. This identity holds for all $\kappa$ and implies the claimed
formula when $\kappa$ is not singular. The formula is again a polynomial in
$\kappa$, hence is valid for all $\kappa.$
\end{proof}
Now we can state the analogue of Corollary 4.5 of \cite{DjO}:
\begin{cor} For $G=G(m,p,N)$,
let $\pi:=x^{({m/p}-1)v}a_\delta$. Let $1$ denote the parameter tuple
such that $1_{C,i}=1$ for both the orbits of reflection hyperplanes $C$, and
$i=1,\dots,e_C-1$ (in the notation
$k=(k_{C,i})$ of Section 2, not the tuple $(\kappa_i)$ of the present section;
cf. Remark \ref{rem:conv}; if $p=m$ or $N=1$ we have only one orbit).
For all $p,q\in P^G$ we have
\begin{equation}
(\pi p,\pi q)_k=(\pi,\pi)_k(p,q)_{k+1}.
\end{equation}
\end{cor}
\begin{proof}
As in \cite{DjO} we have, using Theorem \ref{thm:herm},
Proposition \ref{prop:cycl} and Proposition \ref{prop:sn}
\begin{align*}
(\pi p,\pi q)_k&=(\pi,p^*(T(k))\pi q)_k\\
&=(\pi,\pi p^*(T(k+1))q)_k\\
&=(\pi,\pi)_k(p,q)_{k+1},
\end{align*}
since we may assume that the degrees of $p$ and $q$ are equal.
\end{proof}
This result in fact gives an alternative proof of the explicit description
of the radical which was derived in the previous subsection, analogous to the proof
given in \cite{DjO}.

These results may help in formulating more general statements for arbitrary
complex reflection groups.

\end{document}